\documentclass[reqno,10pt]{amsart}
\usepackage[utf8]{inputenc}
\usepackage[ngerman, english]{babel}
\usepackage[utf8]{inputenc}
\usepackage{pdfsync}
\usepackage{verbatim}
\usepackage[monochrome]{color}
\usepackage{amsmath}
\usepackage{amsthm}
\usepackage{amssymb}
\usepackage{paralist}
\usepackage{pifont}
\usepackage{mathrsfs}
\usepackage{tikz}
\usepackage{tikz-cd}
\usepackage{enumitem}
\numberwithin{equation}{section}
\newtheorem{thm}{Theorem}[section]

  \newtheorem{lemma}[thm]{Lemma}
  \newtheorem{cor}[thm]{Corollary}
  \newtheorem{thmx}{Theorem}
   
  \newtheorem{corx}[thmx]{Corollary}
  \newtheorem{prop}[thm]{Proposition}
  \newtheorem{quest}[thm]{Question}

  \theoremstyle{definition}
 \newtheorem{defi}[thm]{Definition}
  \newtheorem{rem}[thm]{Remark}
    {\begin{proof}[Beweis]}
    {\end{proof}}

\usepackage{tikz}
\usetikzlibrary{matrix,arrows,positioning}
\usepackage{hyperref}
\allowdisplaybreaks
\title[Diagonal and Fine Tower Dimensions on Circles and Products]{Diagonal and Fine Tower Dimensions for Minimal Circle Homeomorphisms and Product Actions}
\author[K.\ Li]{Kang Li}
\address[Kang Li]{Centre for Mathematical Sciences, Lund University, Box 118, SE-221 00 Lund, Sweden}
\email{kang.li@math.lth.se}
\author[G.\ Wang]{Geng Wang}
\address[Geng Wang]{School of Mathematical Sciences, East China Normal University, 500 Dongchuan Road, Minhang District, Shanghai 200241, P.R. China}
\email{51265500012@stu.ecnu.edu.cn}
\date{}


\usepackage{graphicx}

\begin{document}

\begin{abstract}
We study minimal circle homeomorphisms and their products with Cantor systems, and show that these actions have diagonal dimension at most two and fine tower dimension at least two. Consequently, diagonal dimension, fine tower dimension, and tower dimension are all equal to two for minimal circle homeomorphisms. We also obtain a strict separation between tower dimension and fine tower dimension on a one-dimensional space.
\end{abstract}

\maketitle

\section{Introduction}

Nuclear dimension, introduced by Winter and Zacharias \cite{MR2609012}, is a noncommutative analogue of covering dimension. On the dynamical side, Kerr \cite{MR4167017} introduced the tower dimension and the fine tower dimension of an action, which relate open tower decompositions of the underlying space to finite-dimensional approximations of the crossed product. Tower dimension can thus be regarded as a partial dynamical analogue of nuclear dimension. Fine tower dimension strengthens this notion by requiring the tower levels to refine any prescribed finite open cover, and therefore also reflects the covering dimension of the underlying space. 

Li, Liao, and Winter \cite{MR5147990} introduced diagonal dimension, a refinement of nuclear dimension that keeps track of a distinguished abelian subalgebra and its normalizers. For an action $G \curvearrowright X$ of a countable discrete group on a compact Hausdorff space, they showed that 
\[
\dim_{\mathrm{ftow}}(X,G) \le \dim_{\mathrm{diag}}(C(X) \subseteq C(X) \rtimes_{\mathrm r} G),
\]
with equality when $X$ is zero-dimensional. Moreover, finiteness of either dimension forces the action to be free; see \cite[Remark~5.3(iv) and Theorem~5.4]{MR5147990}. It is therefore natural to ask when equality holds beyond the zero-dimensional case.

\begin{quest} \label{ques:A}
	Let $\tau$ be a homeomorphism of a compact Hausdorff space $X$ that generates a free $\mathbb Z$-action. Under what conditions does
	\[
	\dim_{\mathrm{ftow}}(X,\mathbb Z) = \dim_{\mathrm{diag}}(C(X) \subseteq C(X) \rtimes_\tau \mathbb Z)
	\]
	hold?
\end{quest}

Irrational rotations provide a natural test case. Their crossed products are the irrational rotation algebras, which are fundamental examples in $\mathrm{C}^*$-algebra theory; see Rieffel \cite{MR623572} and Connes \cite{MR1303779}. Moreover, every minimal homeomorphism of the circle $\mathbb T = \mathbb R/\mathbb Z$ is topologically conjugate to an irrational rotation; see \cite[Proposition~11.1.4 and Theorem~11.2.7(1)]{MR1326374}. Since the dimensions considered here are invariant under topological conjugacy (see Proposition~\ref{prop:topological-conjugacy}), it suffices to consider irrational rotations when studying minimal circle homeomorphisms. We state our first main result in the more general setting of product actions with an arbitrary compact zero-dimensional factor. 

\begin{thmx} [see Theorem~\ref{thm:upper-bound-diagonal-dimension}] \label{thm:A}
	Let $\tau\colon \mathbb T\to \mathbb T$ be a minimal homeomorphism, and let $\sigma\colon K\to K$ be a homeomorphism of a nonempty compact zero-dimensional Hausdorff space $K$. Equip $\mathbb T\times K$ with the diagonal $\mathbb Z$-action generated by $\tau\times\sigma$. Then
	\[
	\dim_{\mathrm{diag}}(C(\mathbb T\times K)\subseteq C(\mathbb T\times K)\rtimes_{\tau\times\sigma}\mathbb Z)\le 2.
	\]
\end{thmx}

To prove Theorem~\ref{thm:A}, we reduce to an irrational rotation and use the Rokhlin towers arising from the continued-fraction expansion of its rotation number, together with a suitably translated copy of these towers. We then refine the tower bases and select pairwise disjoint levels of arbitrarily small length that cover the circle up to finitely many points (Lemma~\ref{lem:castles-refinement}). An additional tower covers the remaining points (Lemma~\ref{lem:boundary-tower}), and suitable bump functions associated with these towers form a partition of unity. Finite clopen partitions of the zero-dimensional factor extend the construction to the product without adding colors. The outgoing maps are constructed from the towers, following Liao and Tikuisis \cite{MR4339011}, and send matrix units to normalizers. For the incoming maps, we compress regular representations to finite matrices indexed by the tower levels and multiply on both sides by the square roots of the corresponding averaged weights, preserving complete positivity. The incoming maps send the unit to the averaged weights, whose images under the outgoing maps still form a partition of unity. The composite map is therefore unital and completely positive, hence contractive. The partition-of-unity identity also allows us to combine the estimates for the three colors. The approximation error is controlled by the oscillation of the coefficients on the product tower levels and the variation of the averaged weights under fixed integer shifts of the level indices (Lemma~\ref{lem:compression-estimate}), both of which can be made arbitrarily small.

For minimal circle homeomorphisms, the lower bound $\dim_{\mathrm{tow}}(\mathbb T,\mathbb Z)\ge 2$ follows from work of Chinen and Yamauchi \cite[Theorem~1.1]{MR5032552} together with \cite[Theorem~5.14]{MR4167017}. Combining this lower bound with Theorem~\ref{thm:A} applied to a one-point space $K$, we obtain the following answer to Question~\ref{ques:A}.

\begin{corx} [see Corollary~\ref{cor:dimensions-minimal-circle-homeo}] \label{cor:A}
	Let $\tau\colon \mathbb{T}\to \mathbb{T}$ be a minimal homeomorphism, and equip $\mathbb{T}$ with the $\mathbb{Z}$-action generated by $\tau$. Then
	\[
	\dim_{\mathrm{tow}}(\mathbb{T},\mathbb{Z}) = \dim_{\mathrm{ftow}}(\mathbb{T},\mathbb{Z}) = \dim_{\mathrm{diag}}(C(\mathbb{T}) \subseteq C(\mathbb{T}) \rtimes_\tau \mathbb{Z}) = 2.
	\]
\end{corx}

Now consider the product actions in Theorem~\ref{thm:A}. If the action on $K$ is free, the argument for the lower bound on the circle does not carry over: every path in $\mathbb T\times K$ has constant $K$-coordinate, and hence cannot join two distinct points of the same orbit. Our second main result gives a matching lower bound for the fine tower dimension by a different argument, which also applies when the circle is replaced by a more general connected space. 

\begin{thmx} [see Theorem~\ref{thm:product-fine-lower-bound}] \label{thm:B}
	Let $Y$ be an infinite, connected, locally connected compact Hausdorff space, and let $\tau\colon Y\to Y$ be a minimal homeomorphism. Let $K$ be a nonempty compact zero-dimensional Hausdorff space, and let $\sigma\colon K\to K$ be a homeomorphism. Equip $Y\times K$ with the diagonal $\mathbb Z$-action generated by $\tau\times\sigma$. Then
	\begin{enumerate}[label=(\roman*)]
		\item $\dim_{\mathrm{ftow}}(Y\times K,\mathbb Z)\ge 2$;
		\item if $\sigma$ has a periodic point, then $\dim_{\mathrm{tow}}(Y\times K,\mathbb Z)\ge 2$.
	\end{enumerate}
\end{thmx}

The proof of Theorem~\ref{thm:B} uses Lemma~\ref{lem:rectangular-fine-towers} to refine the towers into ones with rectangular bases and connected first factors. In the fine tower dimension argument, the prescribed open cover ensures that these first factors are proper subsets of $Y$. In the tower dimension argument, properness follows from the disjointness of tower levels after restriction to a periodic orbit of $\sigma$ in $K$. In either case, connectedness yields a chain of inclusions among translates of the first factors, leading to a contradiction with total minimality of $\tau$.

As a further application, we compare tower dimension with fine tower dimension. Kerr \cite[Example~9.4]{MR4167017} already exhibited a strict separation between them on an infinite-dimensional space. In finite dimensions, if a $\mathbb Z$-action on a compact Hausdorff space $X$ with $\dim (X) = d < \infty$ admits a free Cantor system as a factor, then the estimates \eqref{eq:fine-diagonal}, \eqref{eq:2.2}, and \eqref{eq:tower-dimension-one} recalled in Section~\ref{sec:preliminaries} yield
\begin{equation} \label{eq:free-cantor-factor-bounds}
	\dim_{\mathrm{tow}}(X,\mathbb Z) = 1\quad \text{and}\quad d\leq \dim_{\mathrm{ftow}}(X,\mathbb Z) \leq 2d+1;
\end{equation}
see \cite[Theorem 2.2]{MR4277767} and \cite[Propositions~4.6, 4.8, and 4.11; Corollary~5.15]{MR4167017}. Hence the two dimensions differ whenever $d\ge 2$, whereas these estimates leave the case $d=1$ open.

The following corollary determines all three dimensions for products of minimal circle homeomorphisms with Cantor systems. In this one-dimensional setting, the tower dimension is strictly smaller than both the fine tower dimension and the diagonal dimension precisely when the Cantor action is free.

\begin{corx} [see Corollary~\ref{cor:products-dimensions}] \label{cor:B}
	Let $\tau\colon \mathbb T \to \mathbb T$ be a minimal homeomorphism, and let $\sigma\colon K \to K$ be a homeomorphism of a Cantor space $K$. Equip $\mathbb T\times K$ with the diagonal $\mathbb Z$-action generated by $\tau\times\sigma$. Then
	\[
	\dim_{\mathrm{ftow}}(\mathbb T\times K,\mathbb Z) = \dim_{\mathrm{diag}}(C(\mathbb{T}\times K) \subseteq C(\mathbb{T}\times K) \rtimes_{\tau\times\sigma} \mathbb{Z}) = 2.
	\]
	Moreover,
	\[
	\dim_{\mathrm{tow}}(\mathbb{T} \times K, \mathbb{Z})=
	\begin{cases}
		1, & \text{if } \sigma \text{ generates a free } \mathbb{Z}\text{-action on } K,\\
		2, & \text{otherwise}.
	\end{cases}
	\]
\end{corx}

The circle case in Corollary~\ref{cor:A} also suggests a question beyond minimality. 

\begin{quest} \label{ques:B}
	Let $\tau\colon \mathbb T\to\mathbb T$ be an orientation-preserving homeomorphism with irrational rotation number, so that the induced $\mathbb Z$-action is free but not necessarily minimal. Is it true that
	\[
	\dim_{\mathrm{ftow}}(\mathbb T,\mathbb Z) = \dim_{\mathrm{diag}}(C(\mathbb T)\subseteq C(\mathbb T)\rtimes_{\tau}\mathbb Z)=2?
	\]
\end{quest}

The lower bound $\dim_{\mathrm{tow}}(\mathbb T,\mathbb Z)\ge 2$ still holds for these actions: the argument in the proof of Corollary~\ref{cor:dimensions-minimal-circle-homeo} only uses that $\mathbb T$ is path-connected and that the action is free. Thus Question~\ref{ques:B} amounts to showing that the diagonal dimension is at most two. Every orientation-preserving circle homeomorphism with irrational rotation number is semiconjugate to an irrational rotation; see \cite[Proposition~6.17]{MR1876932}. However, this semiconjugacy may collapse nontrivial intervals, so our tower construction does not extend directly to the non-minimal case. For crossed products associated with non-minimal circle homeomorphisms, see Putnam, Schmidt, and Skau \cite{MR847334} and Bell \cite{arXiv:2604.18410}.

The paper is organized as follows. In Section~\ref{sec:preliminaries}, we recall the dimension invariants, the relations between them, and their invariance under topological conjugacy. In Section~\ref{sec:diagonal dimension for minimal circle homeomorphisms}, we construct refined tower decompositions of the circle and the associated diagonal approximations, and prove Theorem~\ref{thm:A} and Corollary~\ref{cor:A}. In Section~\ref{sec:fine tower and diagonal dimensions of products with Cantor systems}, we establish a rectangular refinement of tower covers, and prove Theorem~\ref{thm:B} and Corollary~\ref{cor:B}.

\medskip

\noindent{\it Declaration of AI use.} During the development and preparation of this work, the authors used OpenAI's ChatGPT and Anthropic's Claude as exploratory tools. In particular, GPT-5.5 Thinking suggested replacing our initial approach to the upper bound for the fine tower dimension, which was based on orbit disconnection, with a more direct approach using Rokhlin partitions of the circle arising from continued fractions. Both GPT-6 Sol and Claude Fable 5.1 were used to explore a partition of unity built from these partitions in the argument for the diagonal dimension. All references were checked against the original sources. The authors developed and refined the resulting constructions, estimates, and proofs, determined the structure and presentation of the paper, and take full responsibility for its content.

\section{Preliminaries} \label{sec:preliminaries}

Throughout this section, $G$ denotes a countable discrete group acting on a compact Hausdorff space $X$. By a Cantor space we mean a space homeomorphic to the standard Cantor set. 

Let $C_c(G,C(X))$ denote the \textit{algebraic crossed product} of the action $G\curvearrowright X$, whose elements are finite formal sums $\sum_{g\in F} f_gu_g$ with $F\subseteq G$ finite and $f_g\in C(X)$ for each $g\in F$. Its $*$-algebra structure is determined by the usual operations on $C(X)$, the relations $u_gu_h=u_{gh}$, $u_g^*=u_{g^{-1}}$ for $g,h\in G$, and the covariance relation
\[
u_gfu_g^*=f\circ g^{-1},\qquad g\in G,\ f\in C(X),
\]
where $(f\circ g^{-1})(x)=f(g^{-1}x)$. The \textit{full crossed product} $C(X)\rtimes G$ and the \textit{reduced crossed product} $C(X)\rtimes_{\mathrm r} G$ are the completions of $C_c(G,C(X))$ with respect to the maximal and the reduced $\mathrm{C}^*$-norm, respectively. The maximal norm is defined using all covariant representations, while the reduced norm is obtained from the regular representation induced by a faithful representation of $C(X)$. For these constructions and their basic properties, see \cite[Section~4.1]{MR2391387}. We regard $C(X)$ as a unital $\mathrm{C}^*$-subalgebra of each crossed product via the canonical embedding.

Since $\mathbb Z$ is amenable, the full and reduced crossed products of a $\mathbb Z$-action coincide. When the action is generated by a homeomorphism $\tau\colon X\to X$, we write $C(X)\rtimes_\tau \mathbb Z$ for the associated crossed product, omitting the subscript $\mathrm r$.

\subsection{Tower and diagonal dimensions} \label{subsec:tower and diagonal dimensions}

We begin by recalling the tower dimension and the fine tower dimension of a topological dynamical system.

\begin{defi} [\protect{\cite[Definitions~4.1 and 4.2]{MR4167017}}]
	A \textit{tower} is a pair $(V,S)$, where $V\subseteq X$ and $S\subseteq G$ is finite, such that the sets $sV$, $s\in S$, are pairwise disjoint. We call $V$ the \textit{base}, $S$ the \textit{shape}, the sets $sV$ the \textit{levels}, and $SV=\bigcup_{s\in S} sV$ the \textit{support} of the tower. A tower is \textit{open} if its base is open. A family of towers is a \textit{castle} if all of its levels are pairwise disjoint. A family of towers \textit{covers} $X$ if the union of its supports is $X$.
	
	Let $E\subseteq G$ be finite. A family of towers $\{(V_i,S_i)\}_{i\in I}$ covering $X$ is \textit{$E$-Lebesgue} if for every $x\in X$ there exist $i\in I$ and $t\in S_i$ such that $x\in tV_i$ and $Et\subseteq S_i$. We refer to this property as $E$-Lebesgueness.
\end{defi}

For finite subsets $E,S\subseteq G$, the \textit{$E$-interior} of $S$ is 
\[
S_E^\circ = \{s\in S : Es\subseteq S\}.
\]
For a tower $(V,S)$, a level $sV$ is called an \textit{$E$-interior level} if $s\in S_E^\circ$.

\begin{defi} [\protect{\cite[Definition~4.3]{MR4167017} and \cite[Definition~5.2]{MR5147990}}] \label{def:tower-dimension}
	The action $G\curvearrowright X$ is said to have \textit{tower dimension} at most $d$, written $\dim_{\mathrm{tow}}(X,G)\le d$, if for every finite set $E\subseteq G$, there exists an $E$-Lebesgue family of open towers $\{(V_i,S_i)\}_{i\in I}$ whose supports have chromatic number at most $d+1$, that is, the towers can be partitioned into at most $d+1$ castles. Note that we may always take the index set $I$ to be finite in view of the compactness of $X$. 
	
	If, in addition, for any given open cover $\mathcal U$ of $X$ such a family can be chosen so that each level is contained in some member of $\mathcal U$, then the action is said to have \textit{fine tower dimension} at most $d$, written $\dim_{\mathrm{ftow}}(X,G)\le d$. If no such $d$ exists, then the tower dimension (or the fine tower dimension, respectively) is said to be infinite. 
\end{defi}

For metrizable $X$, this open-cover formulation of fine tower dimension agrees with the metric one, in which the levels are required to have arbitrarily small diameter. We do not assume freeness in general, although finite tower dimension implies it; see \cite[Remarks~5.3(ii) and~(iv)]{MR5147990}. 

The two tower dimensions above are defined purely in dynamical terms. We next recall the diagonal dimension of an inclusion of $\mathrm{C}^*$-algebras. For a $\mathrm{C}^*$-algebra $A$, we write $A_+$ for its positive cone. Given an inclusion $D \subseteq A$ of $\mathrm{C}^*$-algebras, the set of \textit{normalizers} of $D$ in $A$ is 
\[
\mathcal{N}_A(D) = \{ a \in A : aDa^* + a^*Da \subseteq D \}.
\]
If $F$ is a finite-dimensional $\mathrm{C}^*$-algebra and $D_F$ is a maximal abelian $*$-subalgebra of $F$, an element $v\in F$ is called a \textit{matrix unit with respect to $D_F$} if $v^*v$ and $vv^*$ are minimal projections in $D_F$. A completely positive map $\varphi\colon A\to B$ between $\mathrm{C}^*$-algebras has \textit{order zero} if $\varphi(a)\varphi(b)=0$ whenever $a,b\in A_+$ satisfy $ab=0$; see \cite{MR2545617} for the structure theory of such maps.

\begin{defi}[\protect{\cite[Definition~2.1]{MR5147990}}] \label{def:diagonal-dimension}
	Let $D_A\subseteq A$ be an inclusion of $\mathrm{C}^*$-algebras with $D_A$ abelian. We say that $D_A\subseteq A$ has \textit{diagonal dimension} at most $d$, written $\dim_{\mathrm{diag}}(D_A\subseteq A)\leq d$, if for every finite subset $\mathfrak{F}\subseteq A$ and every $\varepsilon>0$ there exist a finite-dimensional $\mathrm{C}^*$-algebra $F$ with a maximal abelian $*$-subalgebra $D_F$ and completely positive maps
	\[
	A\xrightarrow{\psi}F\xrightarrow{\varphi}A
	\]
	such that
	\begin{enumerate}[label=(\roman*)]
		\item $\psi$ is contractive;
		\item $\|\varphi\psi(a)-a\|<\varepsilon$ for every $a\in\mathfrak{F}$;
		\item $F$ admits a decomposition $F=\bigoplus_{i=0}^d F^{(i)}$ such that $\varphi^{(i)}:=\varphi|_{F^{(i)}}$ is completely positive contractive (c.p.c.\ for short) order zero for each $i=0,\dots,d$;
		\item $\psi(D_A)\subseteq D_F$;
		\item $\varphi(v)\in \mathcal{N}_A(D_A)$ for every matrix unit $v$ with respect to $D_F$.
	\end{enumerate}
	As usual, we write $\dim_{\mathrm{diag}}(D_A\subseteq A)=d$ if $d$ is the least integer such that $\dim_{\mathrm{diag}}(D_A\subseteq A)\leq d$. If no such $d$ exists, we write $\dim_{\mathrm{diag}}(D_A\subseteq A)=\infty$.
\end{defi}

We refer to the summands $F^{(i)}$ as the colors of the approximation. 

%Recall that an inclusion $D\subseteq A$ of $\mathrm{C}^*$-algebras is said to be \textit{nondegenerate} if $D$ contains an approximate unit for $A$. For an action $G \curvearrowright X$, we consider the canonical inclusion of $C(X)$ into the associated crossed product $C(X) \rtimes_{\mathrm r} G$, which is nondegenerate because $C(X)$ contains the unit of the crossed product.

\subsection{Dimension relations and topological conjugacy} \label{subsec:dimension relations and topological conjugacy}

We next recall the dimension estimates used below. We write $\dim(X)$ for the covering dimension of $X$ and, for each notion of dimension in this paper, use the superscript $+1$ to denote its value plus one; for example, $\dim^{+1}(X)=\dim(X)+1$. 

The estimates in \cite[Proposition~4.11]{MR4167017} and \cite[Theorem~5.4]{MR5147990} give 
\begin{align}
	\max\{\dim_{\mathrm{tow}}^{+1}(X,G),\dim^{+1}(X)\}&\leq \dim_{\mathrm{ftow}}^{+1}(X,G) \notag \\
	&\leq \dim_{\mathrm{diag}}^{+1}(C(X)\subseteq C(X)\rtimes_\mathrm r G) \label{eq:fine-diagonal}\\
	&\leq\dim_{\mathrm{tow}}^{+1}(X,G)\cdot \dim^{+1}(X). \label{eq:2.2}
\end{align}
In particular, when $X$ is zero-dimensional, one has 
\[
\dim_{\mathrm{tow}}(X,G)=\dim_{\mathrm{ftow}}(X,G)=\dim_{\mathrm{diag}}(C(X)\subseteq C(X)\rtimes_\mathrm r G).
\]

We also recall that the tower dimension is monotone under extensions. Suppose that $G\curvearrowright Y$ is a free action on a compact Hausdorff space that is an extension of $G\curvearrowright X$, that is, there is a $G$-equivariant continuous surjection $\pi\colon Y\to X$. Then 
\begin{equation} \label{eq:extension-monotonicity}
	\dim_{\mathrm{tow}}(Y,G)\le \dim_{\mathrm{tow}}(X,G); 
\end{equation}
see \cite[Proposition~4.6]{MR4167017}. 

This applies in particular to $\mathbb Z$-actions with a free Cantor factor. Recall that every free $\mathbb Z$-action on a Cantor space $K$ satisfies $\dim_{\mathrm{tow}}(K,\mathbb Z)=1$ by \cite[Theorem~2.2]{MR4277767} and \cite[Corollary~5.15]{MR4167017}. Consequently, if a $\mathbb Z$-action on a compact Hausdorff space $Y$ admits a free Cantor factor $\mathbb Z \curvearrowright K$, then \eqref{eq:extension-monotonicity} gives $\dim_{\mathrm{tow}}(Y,\mathbb Z)\le 1$. Since $\mathbb Z$ is not locally finite, every $\mathbb Z$-action has tower dimension at least one, and hence 
\begin{equation} \label{eq:tower-dimension-one}
	\dim_{\mathrm{tow}}(Y,\mathbb Z) = 1. 
\end{equation}

Recall that two actions $G\curvearrowright X$ and $G\curvearrowright Y$ are \textit{topologically conjugate} if there is a homeomorphism $\chi\colon X\to Y$ such that $\chi(gx)=g\chi(x)$ for all $g\in G$ and $x\in X$. We next establish the invariance of these three dimensions under topological conjugacy. 

\begin{prop} \label{prop:topological-conjugacy}
	Let $G\curvearrowright X$ and $G\curvearrowright Y$ be topologically conjugate free actions on compact Hausdorff spaces. Then
	\[
	\dim_{\mathrm{tow}}(X,G)=\dim_{\mathrm{tow}}(Y,G) \quad \text{and} \quad \dim_{\mathrm{ftow}}(X,G)=\dim_{\mathrm{ftow}}(Y,G).
	\]
	Moreover, $\dim_{\mathrm{diag}}(C(X)\subseteq C(X)\rtimes_\mathrm r G) = \dim_{\mathrm{diag}}(C(Y)\subseteq C(Y)\rtimes_\mathrm r G)$. 
\end{prop}

\begin{proof}
	Since the actions are topologically conjugate, there is a $G$-equivariant homeomorphism $\chi\colon X\to Y$. We first consider tower dimension. Let $E \subseteq G$ be finite, and let $\{(V_i, S_i)\}_{i\in I}$ be a finite $E$-Lebesgue family of open towers on $Y$. For every $g \in G$ and $i \in I$, we have 
	\[
	g\chi^{-1}(V_i) = \chi^{-1}(gV_i).
	\]
	Hence $\{(\chi^{-1}(V_i), S_i)\}_{i\in I}$ is a family of open towers on $X$, and it is again $E$-Lebesgue. Since two tower supports intersect if and only if their preimages under $\chi$ intersect, the chromatic number is unchanged. It follows that 
	\[
	\dim_{\mathrm{tow}}(X,G)\le \dim_{\mathrm{tow}}(Y,G).
	\]
	For fine tower dimension, given a finite open cover $\mathcal U$ of $X$, choose such a tower family on $Y$ whose levels refine the finite open cover $\{\chi(U) : U \in \mathcal U\}$ of $Y$. Then, for each $i\in I$ and $s\in S_i$, the pullback level satisfies $s\chi^{-1}(V_i) = \chi^{-1}(sV_i) \subseteq U$ for some $U \in \mathcal U$. It follows that 
	\[
	\dim_{\mathrm{ftow}}(X,G)\le \dim_{\mathrm{ftow}}(Y,G).
	\]
	Applying the same arguments to $\chi^{-1}$ yields the reverse inequalities for both dimensions.
	
	Finally, the pullback $\chi^* \colon C(Y) \to C(X)$, defined by $\chi^*(f) = f \circ \chi$, is a $G$-equivariant $*$-isomorphism. It induces a $*$-isomorphism $\widetilde{\chi}\colon C(Y) \rtimes_\mathrm r G \to C(X) \rtimes_\mathrm r G$ satisfying $\widetilde{\chi}(C(Y)) = C(X)$. Suppose that $\dim_{\mathrm{diag}}(C(Y)\subseteq C(Y)\rtimes_\mathrm r G)\le d<\infty$. Given a finite subset $\mathfrak F_X\subseteq C(X)\rtimes_\mathrm r G$ and $\varepsilon>0$, choose a diagonal approximation 
	\[
	C(Y)\rtimes_\mathrm r G \xrightarrow{\psi} F \xrightarrow{\varphi} C(Y)\rtimes_\mathrm r G
	\]
	as in Definition~\ref{def:diagonal-dimension} with $d+1$ colors for $\widetilde{\chi}^{-1}(\mathfrak F_X)$ within $\varepsilon$. Define 
	\[
	\psi_X=\psi\circ \widetilde{\chi}^{-1}\quad\text{and}\quad \varphi_X=\widetilde{\chi}\circ \varphi.
	\]
	These maps are completely positive, $\psi_X$ is contractive, and each restriction of $\varphi_X$ to a color is c.p.c. order zero. Moreover, $\psi_X(C(X))\subseteq D_F$, and $\varphi_X$ sends every matrix unit with respect to $D_F$ to a normalizer of $C(X)$, since $\widetilde{\chi}$ maps normalizers of $C(Y)$ onto normalizers of $C(X)$. As $\widetilde{\chi}$ is isometric, we have 
	\[
	\|\varphi_X\psi_X(a)-a\|=\|\varphi\psi(\widetilde{\chi}^{-1}(a))-\widetilde{\chi}^{-1}(a)\|<\varepsilon,\qquad a\in \mathfrak F_X.
	\]
	Thus $\dim_{\mathrm{diag}}(C(X)\subseteq C(X)\rtimes_\mathrm r G)\le d$. Applying the same argument to $\widetilde{\chi}^{-1}$ gives equality.
\end{proof}

%-------------------------------------------%

\section{Diagonal dimension for minimal circle homeomorphisms} \label{sec:diagonal dimension for minimal circle homeomorphisms}

In this section, we refine the Rokhlin tower decompositions of an irrational rotation given by its continued-fraction expansion, and use them to construct diagonal approximations with three colors for its products with zero-dimensional systems.

Throughout this section, we equip $\mathbb{T} = \mathbb{R}/\mathbb{Z}$ with its standard metric and fix an irrational number $\alpha\in(0,1)$. Let 
\[
R_\alpha(x) = x+\alpha \pmod{1},\qquad x\in\mathbb{T},
\]
be the corresponding irrational rotation. We identify $\mathbb{T}$ with $[0,1)$ in the usual way. 

\subsection{Continued-fraction Rokhlin towers and their refinements} \label{subsec:continued-fraction Rokhlin towers and their refinements}

We denote by $(p_n/q_n)_{n\ge 0}$ the continued-fraction convergents of $\alpha$, indexed so that $p_0/q_0 = 0/1$, and set $\varepsilon_n = q_n \alpha - p_n$. Then $(-1)^n \varepsilon_n > 0$ and $|\varepsilon_n| \to 0$. For every even integer $n \ge 2$, put $a = \varepsilon_n$ and $b = -\varepsilon_{n-1}$. The identity $p_nq_{n-1}-p_{n-1}q_n=(-1)^{n-1}$ yields 
\[
q_{n-1}a + q_nb = 1.
\]
Define the half-open intervals
\begin{equation} \label{eq:I_n,J_n}
	I_n = [0,a) \quad \text{and} \quad J_n = [1-b,1).
\end{equation}
The endpoints of $I_n$ are $0$ and $R_\alpha^{q_n}(0)$, and those of $J_n$ are $R_\alpha^{q_{n-1}}(0)$ and $0$. The following standard partition into two towers comes from consecutive convergents. 

\begin{prop} [\protect{\cite[Lecture~9, Theorem~1]{MR1258087}}] \label{prop:continued‑fraction-rokhlin}
	Consider the irrational rotation $R_\alpha$, the sequence $(p_n/q_n)_{n\ge 0}$ of continued-fraction convergents of $\alpha$, and the half-open intervals $I_n,J_n$ defined in \eqref{eq:I_n,J_n}. For every even integer $n \ge 2$, the family
	\[
	\mathcal{P}_n = \{R_\alpha^j(I_n): 0\le j < q_{n-1}\} \cup \{R_\alpha^j(J_n): 0\le j < q_n\}
	\]
	is a partition of $\mathbb{T}$ into half-open intervals. Moreover, every endpoint of an atom of $\mathcal{P}_n$ lies in the orbit $\operatorname{Orb}_{R_\alpha}(0)$.
\end{prop}

For each even integer $n\ge 2$, we pass from the half-open partition $\mathcal P_n$ to open towers by taking interiors of the base intervals. Set 
\[
U_{n,0} = \operatorname{int}(I_n)\quad\text{and}\quad U_{n,1} = \operatorname{int}(J_n),
\]
and define 
\[
S_{n,0} = \{0,1,\dots,q_{n-1}-1\}\quad \text{and}\quad S_{n,1} = \{0,1,\dots,q_n-1\}.
\]
Write $D_n$ for the finite set of endpoints of the intervals in $\mathcal{P}_n$, namely 
\begin{equation} \label{notation: D_n}
	D_n = \{R_\alpha^k(0): 0\le k< q_n+q_{n-1}\} \subseteq \operatorname{Orb}_{R_\alpha}(0).
\end{equation}
By Proposition~\ref{prop:continued‑fraction-rokhlin}, the family 
\begin{equation} \label{notation: open castle (U_n0,S_n0)}
	\{(U_{n,0},S_{n,0}),(U_{n,1},S_{n,1})\}
\end{equation}
is an open castle whose levels partition $\mathbb T\setminus D_n$. Since $R_\alpha$ is an isometry, every interval in $\mathcal P_n$ has diameter at most $\max\{a,b\}$. Along the even integers, this bound tends to zero and both tower heights $q_{n-1}$ and $q_n$ tend to infinity.

In the following tower constructions, $E\subseteq \mathbb Z$ will be a finite symmetric subset containing $0$. The $E$-interior levels of $\{(U_{n,0},S_{n,0}),(U_{n,1},S_{n,1})\}$ do not cover $\mathbb T\setminus D_n$ in general. The next lemma uses a translated copy of this castle so that the $E$-interior levels of two castles cover $\mathbb T$ up to finitely many points.

\begin{lemma} \label{lem:translated-castles}
	Consider the endpoint set $D_n$ defined in \eqref{notation: D_n} and the open castle $\{(U_{n,0},S_{n,0}),(U_{n,1},S_{n,1})\}$ defined in \eqref{notation: open castle (U_n0,S_n0)}. Let $E\subseteq\mathbb{Z}$ be a finite symmetric subset containing $0$, and put $m=\max_{e\in E} e$. Choose an integer $N\ge 2m+1$ and an even integer $n\ge 2$ such that both towers have more than $N+2m$ levels. Set
	\[
	W=\bigcup_{t=0}^{N}R_\alpha^t(D_n).
	\]
	Consider the two open castles
	\[
	\mathscr{C}_n^{(0)}=\{(U_{n,0},S_{n,0}),(U_{n,1},S_{n,1})\} \quad\text{and}\quad \mathscr{C}_n^{(1)}=\{(R_\alpha^N U_{n,0},S_{n,0}),(R_\alpha^N U_{n,1},S_{n,1})\}.
	\]
	Then the $E$-interior levels of the towers in $\mathscr{C}_n^{(0)}\cup \mathscr{C}_n^{(1)}$ cover $\mathbb T\setminus W$.
\end{lemma}

\begin{proof}
	For convenience, write
	\[
	S_{n,i}=\{0,1,\dots,s_{n,i}\},\qquad i=0,1.
	\]
	By the choice of $n$, we have $s_{n,i}\ge N+2m$. Since $E\subseteq [-m,m]$, an index $r\in S_{n,i}$ satisfies $E+r\subseteq S_{n,i}$ whenever $m\le r\le s_{n,i}-m$. 
	
	Let $x\in\mathbb{T}\setminus W$. Then $R_\alpha^{-t}(x)\notin D_n$ for $0\le t\le N$. In particular, there are unique $i\in\{0,1\}$ and $r\in S_{n,i}$ such that $x\in R_\alpha^r(U_{n,i})$. We distinguish three cases. 
	
	If $m\le r \le s_{n,i}-m$, then $E+r\subseteq S_{n,i}$, so $x$ lies in an $E$-interior level of $\mathscr{C}_n^{(0)}$. 
	
	If $s_{n,i}-m+1\le r \le s_{n,i}$, then the choice of $n$ and $N$ gives $m<r-N\le s_{n,i}-m$. Hence $E+(r-N)\subseteq S_{n,i}$. Since $x\in R_\alpha^{r-N}(R_\alpha^N U_{n,i})$, the point $x$ lies in an $E$-interior level of $\mathscr{C}_n^{(1)}$.
	
	Finally, if $0\le r\le m-1$, we show that $x$ lies in an $E$-interior level of $\mathscr{C}_n^{(1)}$. Put $z=R_\alpha^{-r}(x)\in U_{n,i}$. Since $r+1\le m\le N$, we have 
	\[
	R_\alpha^{-1}(z)=R_\alpha^{-(r+1)}(x)\notin D_n.
	\]
	Thus $R_\alpha^{-1}(z)$ lies in some level of $\mathscr{C}_n^{(0)}$. This must be the top level of its tower: otherwise $z$ would also lie in a level of $\mathscr{C}_n^{(0)}$ with positive index, contradicting the disjointness of the levels. Hence $R_\alpha^{-1}(z)\in R_\alpha^{s_{n,j}}(U_{n,j})$ for some $j\in\{0,1\}$, and consequently, 
	\[
	R_\alpha^{-N}(x)=R_\alpha^{-(N-r-1)}(R_\alpha^{-1}(z))\in R_\alpha^{s_{n,j}-(N-r-1)}(U_{n,j}).
	\]
	Set $u=s_{n,j}-(N-r-1)$. By the choice of $n$ and $N$, we have  
	\[
	m\le 2m+r< u\quad \text{and}\quad m\le N-r-1= s_{n,j}-u. 
	\]
	Thus $m\le u\le s_{n,j}-m$, and hence $E+u\subseteq S_{n,j}$. Since $x\in R_\alpha^{u}(R_\alpha^N U_{n,j})$, the point $x$ lies in an $E$-interior level of $\mathscr{C}_n^{(1)}$.
\end{proof}

The $E$-interior levels of the two castles in Lemma~\ref{lem:translated-castles} may overlap. The next lemma cuts their bases into finitely many subintervals, which allows us to select pairwise disjoint $E$-interior levels that still cover $\mathbb T$ up to finitely many points of the orbit $\operatorname{Orb}_{R_\alpha}(0)$. 

\begin{lemma} \label{lem:castles-refinement} 
	Let $E\subseteq \mathbb Z$ be a finite symmetric subset containing $0$, and let $\delta>0$. Choose $N$ and an even integer $n\ge 2$ as in Lemma~\ref{lem:translated-castles} such that, in addition, every level of the two castles has diameter less than $\delta$, and let $\mathscr{C}_n^{(0)},\mathscr{C}_n^{(1)}$, and $W$ denote the corresponding castles and finite exceptional set. Then there exist finite open castles $\mathscr R_0,\mathscr R_1$ and a finite set $P$ with $W\subseteq P\subseteq \operatorname{Orb}_{R_\alpha}(0)$ such that the following hold:
	\begin{enumerate}[label=(\roman*)]
		\item for $i=0,1$, each tower $(R,S)\in\mathscr R_i$ is obtained from a tower $(V,S)\in \mathscr C_n^{(i)}$ by restricting its base to an open subinterval $R\subseteq V$ and keeping its shape $S$;
		\item one can simultaneously select a nonempty family of $E$-interior levels from each tower in $\mathscr{R}_0\cup\mathscr{R}_{1}$ such that all selected levels are pairwise disjoint and their union is $\mathbb T\setminus P$.
	\end{enumerate}
\end{lemma}

We fix such a selection. For each refined tower $(R,S)$, we write $A_{R,S}\subseteq S_E^\circ$ for the nonempty set of indices of its selected levels. We call $A_{R,S}$ the set of \textit{active indices} of $(R,S)$, and the levels $R_\alpha^s(R)$ with $s\in A_{R,S}$ its \textit{active levels}. Thus Lemma~\ref{lem:castles-refinement}(ii) states that 
\[
\mathbb T \setminus P = \bigsqcup_{i=0}^{1}\bigsqcup_{(R,S) \in \mathscr{R}_i}\bigsqcup_{s \in A_{R,S}} R_{\alpha}^{s}(R).
\]

\begin{proof} [Proof of Lemma~\ref{lem:castles-refinement}]
	Let $\mathscr{J}$ be the family of connected components of $\mathbb{T}\setminus W$. By Lemma~\ref{lem:translated-castles}, the $E$-interior levels of the towers in $\mathscr{C}_n^{(0)}\cup \mathscr{C}_n^{(1)}$ cover $\mathbb T\setminus W$, and the boundary of each such level is contained in $W$. Hence each $J\in\mathscr{J}$ is contained in one of these levels. 
	
	For each $J\in\mathscr J$, fix a tower $(V_J,S_J)\in \mathscr{C}_n^{(0)}\cup\mathscr{C}_n^{(1)}$ and an index $s_J\in (S_J)_E^\circ$ such that $J\subseteq R_\alpha^{s_J}(V_J)$. For each original tower $(V,S)$, let 
	\[
	\mathscr{J}(V,S)=\{J\in\mathscr{J} : (V_J,S_J)=(V,S)\}
	\]
	be the family of components assigned to $(V,S)$, and define the finite cutting set
	\[
	\Xi(V,S)=\partial V\cup\bigcup_{J\in\mathscr{J}(V,S)}\partial(R_\alpha^{-s_J}(J)).
	\]
	Note that every connected component of $V\setminus \Xi(V,S)$ avoids $\partial(R_\alpha^{-s_J}(J))$, and is therefore either contained in or disjoint from $R_\alpha^{-s_J}(J)$, for each $J\in\mathscr J(V,S)$. Let $\mathcal{R}(V,S)$ be the family of those components of $V\setminus\Xi(V,S)$ that are contained in $R_\alpha^{-s_J}(J)$ for some $J\in\mathscr{J}(V,S)$. This is a finite family of pairwise disjoint open subintervals of $V$. For $i=0,1$, define
	\[
	\mathscr{R}_i = \{(R,S) : (V,S) \in \mathscr{C}_n^{(i)}\text{ and }R \in \mathcal{R}(V,S)\}.
	\]
	Each tower $(R,S)\in\mathscr R_i$ is obtained from an original tower $(V,S)$ by restricting its base to the open subinterval $R$ and keeping the shape $S$. Within each $\mathscr R_i$, refined levels contained in distinct levels of $\mathscr{C}_n^{(i)}$ are disjoint, and refined levels contained in the same original level are translates of pairwise disjoint refined bases by the same rotation. Hence each $\mathscr R_i$ is a finite open castle, which proves (i).
	
	To specify the selection in (ii), for each $R\in\mathcal{R}(V,S)$, set 
	\[
	A_{R,S}=\{s_J: J\in\mathscr{J}(V,S)\text{ and }R\subseteq R_\alpha^{-s_J}(J)\}.
	\]
	By the definition of $\mathcal{R}(V,S)$ and the choice of $s_J$, each $A_{R,S}$ is a nonempty subset of $S_E^\circ$. We select the levels $R_\alpha^s(R)$ with $s\in A_{R,S}$. For each $J\in\mathscr{J}$, we have $\partial(R_\alpha^{-s_J}(J))\subseteq \Xi(V_J,S_J)$. Consequently, 
	\[
	J\setminus R_\alpha^{s_J}(\Xi(V_J,S_J))=\bigsqcup_{\substack{R\in\mathcal{R}(V_J,S_J)\\ R\subseteq R_\alpha^{-s_J}(J)}}R_\alpha^{s_J}(R).
	\]
	Each selected level $R_\alpha^s(R)$ is contained in a unique component $J\in\mathscr{J}$, namely the one with $(V_J,S_J)=(V,S)$, $s_J=s$, and $R\subseteq R_\alpha^{-s_J}(J)$; we say that this level is \textit{assigned} to $J$. For fixed $J$, the selected levels assigned to $J$ are $R_\alpha^{s_J}$-translates of pairwise disjoint members of $\mathcal{R}(V_J,S_J)$, and hence are pairwise disjoint. Selected levels assigned to different members of $\mathscr J$ are disjoint because these components are disjoint. Thus all selected levels are pairwise disjoint. Define 
	\[
	P=\mathbb T\setminus \bigsqcup_{i=0}^{1}\bigsqcup_{(R,S)\in\mathscr{R}_i}\bigsqcup_{s\in A_{R,S}}R_\alpha^s(R).
	\]
	Since every selected level is contained in $\mathbb T\setminus W$, we have $W\subseteq P$. The decomposition above also shows that $P \subseteq W \cup \bigcup_{J \in \mathscr J} R_{\alpha}^{s_J}(\Xi(V_J, S_J))$, so $P$ is finite. Finally, each point of $\Xi(V,S)$ is either an endpoint of $V$ or an $R_\alpha$-translate of an endpoint of some $J\in\mathscr J(V,S)$. Since $\partial V\subseteq W$, $\partial J\subseteq W$, and $W\subseteq \operatorname{Orb}_{R_\alpha}(0)$, all cutting points and their $R_\alpha$-translates lie in $\operatorname{Orb}_{R_\alpha}(0)$. Hence $W\subseteq P\subseteq \operatorname{Orb}_{R_\alpha}(0)$.
\end{proof}

The following lemma provides an additional tower whose $E$-interior levels cover the remaining finite set $P$ in Lemma~\ref{lem:castles-refinement}. 

\begin{lemma} \label{lem:boundary-tower}
	Let $E\subseteq \mathbb Z$ be a finite symmetric subset containing $0$, let $\delta>0$, and let $P\subseteq \operatorname{Orb}_{R_\alpha}(0)$ be the finite set obtained in Lemma~\ref{lem:castles-refinement} with $E$ and $\delta$. Write 
	\[
	P=\{R_\alpha^\ell(0):\ell\in L\}
	\]
	for a finite set $L\subseteq \mathbb Z$. Set $S_2=L-E$. Then $L\subseteq (S_2)_E^\circ$, and there exists $r>0$ such that $(B_r(0),S_2)$ is an open tower whose levels have diameter less than $\delta$ and whose levels indexed by $L$ cover $P$. Here $B_r(0)$ denotes the open ball of radius $r$ centered at $0$ in $\mathbb T$. 
\end{lemma}

\begin{proof}
	Since $R_\alpha$ generates a free $\mathbb Z$-action, the points $R_\alpha^s(0)$ with $s\in S_2$ are pairwise distinct. Hence we can choose $r>0$ sufficiently small that the sets $R_\alpha^s(B_r(0))$, $s\in S_2$, are pairwise disjoint and have diameter less than $\delta$. Thus $(B_r(0),S_2)$ is an open tower. Since $E=-E$ and $0\in E$, we have $L\subseteq (S_2)_E^\circ$. By the choice of $L$, the levels indexed by $L$ cover $P$.
\end{proof}

We now use the two refined castles $\mathscr R_0,\mathscr R_1$ from Lemma~\ref{lem:castles-refinement} and the tower $(B_r(0),S_2)$ from Lemma~\ref{lem:boundary-tower} to construct a partition of unity of $\mathbb T$.

\begin{lemma} \label{lem:partition-of-unity}
	Let $\mathscr{R}_0,\mathscr{R}_1$ be the refined castles from Lemma~\ref{lem:castles-refinement}, with selected index sets $A_{R,S}$ and the finite set $P=\{R_\alpha^\ell(0):\ell\in L\}$, and let $(B_r(0),S_2)$ be the open tower from Lemma~\ref{lem:boundary-tower} for this $P$. For each base $R$ of a tower in $\mathscr{R}_0\cup\mathscr{R}_1$, define $f_R,\rho\in C(\mathbb T)$ by 
	\[
	f_R(x)=\min\left\{1,\,\frac{d(x,\mathbb{T}\setminus R)}{r}\right\} \quad \text{and} \quad
	\rho(x)=\max\left\{1-\frac{d(x,0)}{r},\,0\right\}.
	\]
	For $t\in\mathbb Z$, let $\theta_t\in\operatorname{Aut}(C(\mathbb T))$ be given by $\theta_t(f)=f\circ R_\alpha^{-t}$, and set 
	\begin{equation*}
		\begin{aligned}
			q_i'&=\sum_{(R,S)\in\mathscr R_i}\sum_{s\in A_{R,S}} \theta_s(f_R),\qquad i=0,1,\\
			q_2'&=\sum_{\ell\in L}\theta_\ell(\rho).
		\end{aligned}
	\end{equation*}
	Then $q_0',q_1',q_2'$ are positive contractions in $C(\mathbb T)$ satisfying
	\[
	q_0'+q_1'+q_2'=1.
	\]
\end{lemma}

\begin{proof}
	By Lemma~\ref{lem:castles-refinement}, the selected levels are pairwise disjoint open intervals whose union is $\mathbb T\setminus P$. Since each of them is connected, open, and disjoint from the others, they are precisely the connected components of $\mathbb T\setminus P$. In particular, the endpoints of each selected level lie in $P$.
	
	Let $x\in \mathbb T\setminus P$, and let $R_\alpha^s(R)$ be the unique selected level containing $x$. All other terms of $q_0'+q_1'$ vanish at $x$, because $\theta_{s'}(f_{R'})$ vanishes outside $R_\alpha^{s'}(R')$. Hence 
	\[
	(q'_0+q'_1)(x) = \theta_s(f_R)(x) = \min\left\{1,\,\frac{d(x,\,\mathbb{T}\setminus R_\alpha^s(R))}{r}\right\} = \min\left\{1,\,\frac{d(x,P)}{r}\right\}.
	\]
	Both sides of this identity vanish on $P$, so it holds on all of $\mathbb T$. 
	
	By Lemma~\ref{lem:boundary-tower}, the balls $R_\alpha^\ell(B_r(0))=B_r(R_\alpha^\ell(0))$, $\ell\in L$, are pairwise disjoint; hence distinct points of $P$ have distance at least $2r$. Suppose that $d(x,P)<r$, and let $\ell_x\in L$ be the unique index with $x\in R_\alpha^{\ell_x}(B_r(0))$. Then $d(x,R_\alpha^\ell(0))>r$ for every $\ell\neq\ell_x$, so $d(x,P)=d(x,R_\alpha^{\ell_x}(0))$ and only the term with $\ell=\ell_x$ contributes to $q_2'(x)$. Therefore, 
	\[
	\begin{aligned}
		q'_2(x) = \theta_{\ell_x}(\rho)(x) 
		= \max\left\{1-\frac{d(x,R_\alpha^{\ell_x}(0))}{r},\,0\right\} 
		= \max\left\{1-\frac{d(x,P)}{r},\,0\right\}.
	\end{aligned}
	\]
	Both sides vanish when $d(x,P)\ge r$, so this identity also holds on all of $\mathbb T$.
	
	Since $\min\{1,t\}+\max\{1-t,0\}=1$ for every $t\ge 0$, the two identities give $q_0'+q_1'+q_2'=1$. As each $q_i'$ is nonnegative, all three are positive contractions.
\end{proof}

\subsection{Diagonal approximations for products with zero-dimensional systems} \label{subsec:diagonal approximations for products with zero-dimensional systems}

In this subsection, we use the refined towers and the partition of unity from the preceding subsection to prove the following theorem.

\begin{thm} \label{thm:upper-bound-diagonal-dimension}
	Let $\tau\colon \mathbb T\to \mathbb T$ be a minimal homeomorphism, and let $\sigma\colon K\to K$ be a homeomorphism of a nonempty compact zero-dimensional Hausdorff space $K$. Equip $\mathbb T\times K$ with the diagonal $\mathbb Z$-action generated by $\tau\times\sigma$. Then
	\[
	\dim_{\mathrm{diag}}(C(\mathbb T\times K)\subseteq C(\mathbb T\times K)\rtimes_{\tau\times\sigma}\mathbb Z)\le 2.
	\]
\end{thm}

Fix $K$ and $\sigma$ as in Theorem~\ref{thm:upper-bound-diagonal-dimension}. We first consider the diagonal $\mathbb Z$-action on $\mathbb T\times K$ generated by $R_\alpha\times\sigma$, and denote by $(u_t)_{t\in\mathbb Z}$ the canonical unitaries in the associated crossed product $C(\mathbb T\times K)\rtimes_{R_\alpha\times\sigma}\mathbb Z$. Thus $u_su_t=u_{s+t}$, $u_t^*=u_{-t}$, and $u_t h u_t^*=h\circ(R_\alpha\times\sigma)^{-t}$ for all $h\in C(\mathbb T\times K)$ and $s,t\in\mathbb Z$.

\medskip
\noindent\textbf{Product towers and partition of unity.} 
Fix a finite symmetric subset $E\subseteq \mathbb Z$ containing $0$, a number $0<\delta<1/4$, and a finite clopen partition $\mathcal P$ of $K$. Let $\mathscr R_0,\mathscr R_1$ be the refined castles from Lemma~\ref{lem:castles-refinement} for $E$ and $\delta$, with selected index sets $A_{R,S}$ and the finite set $P=\{R_\alpha^\ell(0):\ell\in L\}$, and let $(B_r(0),S_2)$ be the additional open tower from Lemma~\ref{lem:boundary-tower} for this $P$ with $S_2=L-E$. Let $f_{R},\rho$ and $q_0',q_1',q_2'$ be the functions associated with these towers in Lemma~\ref{lem:partition-of-unity}. Let 
\[
H= S_2 \cup \bigcup_{i=0}^{1}\bigcup_{(R,S)\in\mathscr R_i} S
\]
be the union of all shapes involved, and choose a finite clopen partition $\mathcal Q$ of $K$ that refines $\sigma^{-s}(\mathcal P)$ for every $s\in H$; equivalently, for all $D\in\mathcal Q$ and $s\in H$, the set $\sigma^s(D)$ is contained in some member of $\mathcal P$. 

Consider the following three families of product towers:
\begin{equation} \label{eq:product-towers}
	\begin{aligned}
		&\{(R\times D,S) : (R,S)\in\mathscr R_i,D\in\mathcal Q\},\qquad i=0,1,\\
		&\{(B_r(0)\times D,S_2) : D\in\mathcal Q\}.
	\end{aligned}
\end{equation}

Each of these three families is an open castle. Indeed, distinct circle levels within each family are disjoint. For each fixed $t\in \mathbb Z$, the family $\{\sigma^t(D):D\in\mathcal Q\}$ is a clopen partition of $K$. Hence the product levels are pairwise disjoint within each family of product towers.

Using the canonical identification $C(\mathbb T\times K)\cong C(\mathbb T)\otimes C(K)$, define 
\begin{equation} \label{eq:f_{R,D}-rho}
	f_{R,D}=f_R\otimes \mathbf 1_D\quad \text{and}\quad \rho_D=\rho\otimes \mathbf 1_D,
\end{equation}
where $\mathbf 1_D\in C(K)$ is the characteristic function of $D$. These are positive contractions with $\operatorname{supp}^\circ(f_{R,D}) \subseteq R\times D$ and $\operatorname{supp}^\circ(\rho_D) \subseteq B_r(0)\times D$, where $\operatorname{supp}^\circ(h)=\{(y,k)\in\mathbb T\times K:h(y,k)\neq 0\}$ denotes the open support of $h$. Set 
\[
q_i=q_i'\otimes \mathbf 1_K,\qquad i=0,1,2.
\]
By Lemma~\ref{lem:partition-of-unity}, these are positive contractions satisfying 
\[
q_0 + q_1 + q_2 = (q_0' + q_1' + q_2') \otimes \mathbf 1_K =1.
\]
With $(\theta_t)_{t\in\mathbb Z}$ as in Lemma~\ref{lem:partition-of-unity}, the covariance relation for the product action gives 
\begin{equation} \label{eq:product-covariance}
	u_t(f\otimes\mathbf 1_D)u_t^*=\theta_t(f)\otimes \mathbf 1_{\sigma^t(D)},\qquad f\in C(\mathbb T),\ D\in \mathcal Q,\ t\in \mathbb Z.
\end{equation}
Since $\sum_{D\in\mathcal{Q}} \mathbf 1_{\sigma^t(D)} = \mathbf 1_K$ for every $t\in\mathbb Z$, we can expand 
\begin{equation} \label{eq:sum-form-q_i}
	\begin{aligned}
		q_i&=\sum_{(R,S)\in\mathscr R_i}\sum_{D\in\mathcal Q}\sum_{s\in A_{R,S}} \theta_s(f_R)\otimes \mathbf 1_{\sigma^s(D)},\qquad i=0,1,\\
		q_2&=\sum_{D\in\mathcal Q}\sum_{\ell\in L} \theta_\ell(\rho)\otimes \mathbf 1_{\sigma^\ell(D)}.
	\end{aligned}
\end{equation}

\medskip
\noindent\textbf{Finite-dimensional algebras and outgoing maps.}
For a single tower, we use the construction of Liao and Tikuisis \cite{MR4339011} to define an order zero map from a matrix algebra into $C(\mathbb T\times K)\rtimes_{R_\alpha\times\sigma}\mathbb Z$.

Let $(V,S)$ be an open tower for the action generated by $R_\alpha\times \sigma$, and let $h\in C(\mathbb{T}\times K)_+$ be a positive contraction with $\operatorname{supp}^\circ(h)\subseteq V$. Set 
\[
M_S=\mathcal{B}(\ell^2(S))\quad \text{and}\quad D_S=\operatorname{span}\{e_{s,s}: s\in S\},
\]
where $(e_{s,t})_{s,t\in S}$ are the standard matrix units of $M_S$, so that $D_S$ is the diagonal subalgebra of $M_S$. Taking all phase functions in \cite[Equation~(3.17)]{MR4339011} to be identically one, we define 
\begin{equation} \label{eq:single-tower-castle-map}
	\Phi_{V,S,h}\colon M_S\to C(\mathbb T\times K)\rtimes_{R_\alpha\times\sigma} \mathbb Z,\qquad \Phi_{V,S,h}(e_{s,t})=u_s h u_t^*,
\end{equation}
and extend linearly. 

We claim that the map $\Phi_{V,S,h}$ is c.p.c.\ order zero and sends each standard matrix unit to a normalizer of $C(\mathbb T\times K)$. The construction in \cite{MR4339011} is stated for compact metrizable spaces, but metrizability is not needed in the single-tower case; for the reader's convenience, we verify these properties directly.

To prove the claim, write $\Phi=\Phi_{V,S,h}$ and $v_s=u_sh^{1/2}$ for $s\in S$. The matrix $[\Phi(e_{s,t})]_{s,t\in S}=[v_sv_t^*]_{s,t\in S}$ is positive, so $\Phi$ is completely positive by Choi's criterion; see \cite[Proposition~1.5.12]{MR2391387}. Since $h$ is a positive contraction and the levels of the tower are pairwise disjoint, we have $\|\Phi\|=\|\sum_{s\in S}\Phi(e_{s,s})\|=\|\sum_{s\in S}u_shu_s^*\|\le 1$, where the first equality follows from the norm formula for positive maps on a unital $\mathrm{C}^*$-algebra; see \cite[Corollary~2.9]{MR1976867}. Moreover, for $t,p\in S$ with $t\neq p$, we have $V\cap (R_\alpha\times\sigma)^{p-t}(V)=\emptyset$, so the covariance relation yields $v_t^*v_p=\delta_{t,p}h$ for all $t,p\in S$. Expanding in matrix units, we obtain $\Phi(x)\Phi(y)=\sum_{s,t\in S}(xy)_{s,t}u_sh^2u_t^*$ for $x,y\in M_S$. In particular, $\Phi(x)\Phi(y)=0$ whenever $x,y\in (M_S)_+$ satisfy $xy=0$, so $\Phi$ has order zero. Finally, both $h$ and the canonical unitaries are normalizers of $C(\mathbb T\times K)$, and hence so is each $\Phi(e_{s,t})=u_shu_t^*$.

We now apply \eqref{eq:single-tower-castle-map} to the three families of product towers in \eqref{eq:product-towers}. Define the finite-dimensional algebras and their diagonal subalgebras by 
\begin{equation} \label{eq:finite-dimensional-algebras}
	\begin{aligned}
		F^{(i)}&=\bigoplus_{(R,S)\in\mathscr R_i}\bigoplus_{D\in\mathcal Q}M_S,
		&
		D_F^{(i)}&=\bigoplus_{(R,S)\in\mathscr R_i}\bigoplus_{D\in\mathcal Q}D_S, \qquad i=0,1,\\
		F^{(2)}&=\bigoplus_{D\in\mathcal Q}M_{S_2},
		&
		D_F^{(2)}&=\bigoplus_{D\in\mathcal Q}D_{S_2}.
	\end{aligned}
\end{equation}
For $i=0,1,2$, each $D_F^{(i)}$ is a maximal abelian $*$-subalgebra of $F^{(i)}$, and we define $\varphi^{(i)}\colon F^{(i)}\to C(\mathbb T\times K)\rtimes_{R_\alpha\times\sigma}\mathbb Z$ by
\begin{equation} \label{eq:color-maps}
	\begin{aligned}
		\varphi^{(i)}(x^{(i)})&=\sum_{(R,S)\in\mathscr R_i}\sum_{D\in\mathcal Q}\Phi_{R\times D,S,f_{R,D}}(x_{R,S,D}), \qquad i=0,1,\\
		\varphi^{(2)}(x^{(2)})&=\sum_{D\in\mathcal Q}\Phi_{B_r(0)\times D,S_2,\rho_D}(x_D),
	\end{aligned}
\end{equation}
where $x^{(i)}=(x_{R,S,D})_{(R,S)\in\mathscr{R}_i,\,D\in\mathcal{Q}}\in F^{(i)}$ for $i=0,1$ and $x^{(2)}=(x_D)_{D\in\mathcal{Q}}\in F^{(2)}$. 

By the claim proved above, each summand in \eqref{eq:color-maps} is c.p.c.\ order zero. Since each family of product towers is an open castle, the maps within each sum have pairwise orthogonal ranges. Hence each $\varphi^{(i)}$ is c.p.c.\ order zero. Moreover, let $v_i\in F^{(i)}$ be a matrix unit with respect to $D_F^{(i)}$. Since $v_i^*v_i$ and $v_iv_i^*$ are minimal projections, exactly one matrix block of $v_i$ is nonzero, and that block has the form $\lambda e_{s,t}$ for some $\lambda\in\mathbb C$ with $|\lambda|=1$. Thus $\varphi^{(i)}(v_i)=\lambda \Phi(e_{s,t})$, where $\Phi$ is the map associated to that matrix block. It follows that 
\begin{equation} \label{eq:matrix-unit-normalizer}
	\varphi^{(i)}(v_i) \in \mathcal N_{C(\mathbb T\times K)\rtimes_{R_\alpha\times\sigma}\mathbb Z}(C(\mathbb T \times K)), \qquad i=0,1,2.
\end{equation}

\medskip
\noindent\textbf{Completely positive incoming maps.}
We next define completely positive maps from the crossed product into the finite-dimensional algebras. To this end, we first construct diagonal weights by averaging translates of the characteristic functions of $A_{R,S}$ in Lemma~\ref{lem:castles-refinement} and $L$ in Lemma~\ref{lem:boundary-tower}.

We use a finite symmetric set $E_0\subseteq \mathbb Z$ containing $0$ to index the Fourier terms in the approximation estimates below. Put $m=\max_{n\in E_0}|n|$, and fix an integer $M\ge 1$ for the averaging. Assume that the set $E$ used in the preceding tower construction satisfies 
\begin{equation} \label{eq:enlargement-E}
	[-M-m,M+m]\cap \mathbb Z\subseteq E.
\end{equation}

Let $\mathbf 1_{A_{R,S}}$ and $\mathbf 1_{L}$ denote the characteristic functions of $A_{R,S}$ and $L$ on $\mathbb Z$. For $s,\ell\in \mathbb Z$, define
\begin{equation} \label{eq:separated-weights}
	\begin{aligned}
		w_{R,S}(s)&=\frac{1}{2M+1}\sum_{j=-M}^{M} \mathbf 1_{A_{R,S}}(s-j),\qquad (R,S)\in\mathscr{R}_0\cup\mathscr{R}_1,\\
		w_{B_r(0),S_2}(\ell)&=\frac{1}{2M+1}\sum_{j=-M}^{M} \mathbf 1_{L}(\ell-j).
	\end{aligned}
\end{equation}
Set
\begin{equation} \label{eq:digonal-weights}
	\begin{aligned}
		W_{R,S}&=\sum_{s\in S} w_{R,S}(s)e_{s,s}\in D_S,\qquad (R,S)\in\mathscr{R}_0\cup\mathscr{R}_1,\\
		W_{B_r(0),S_2}&=\sum_{\ell\in S_2} w_{B_r(0),S_2}(\ell)e_{\ell,\ell}\in D_{S_2}.
	\end{aligned}
\end{equation}
These are diagonal positive contractions, since all their diagonal entries lie in $[0,1]$. 

Since $A_{R,S}\subseteq S_E^\circ$ and $L\subseteq (S_2)_E^\circ$, \eqref{eq:enlargement-E} gives 
\begin{equation} \label{eq:weight-supports}
	\begin{aligned}
		\{s\in\mathbb Z : w_{R,S}(s)>0\}&=A_{R,S}+([-M,M]\cap\mathbb Z)\subseteq S,\\
		\{\ell\in\mathbb Z : w_{B_r(0),S_2}(\ell)>0\}&=L+([-M,M]\cap\mathbb Z)\subseteq S_2.
	\end{aligned}
\end{equation}
Moreover, fix $n\in E_0$. If $w_{R,S}(s)>0$, write $s=t+j$ with $t\in A_{R,S}$ and $-M\le j\le M$. Then $j-n\in E$, so $s-n=t+(j-n)\in A_{R,S}+E\subseteq S$. The same argument applies with $L$ in place of $A_{R,S}$. Thus, for every $n\in E_0$, we have 
\begin{equation} \label{eq:weight-columns}
	\begin{aligned}
		s-n\in S&\quad\text{if}\ w_{R,S}(s)>0,\\
		\ell-n\in S_2&\quad\text{if}\ w_{B_r(0),S_2}(\ell)>0.
	\end{aligned}
\end{equation}
For $s,\ell\in\mathbb Z$, reindexing the defining sums gives 
\begin{equation*}
	\begin{aligned}
		&|w_{R,S}(s)-w_{R,S}(s-1)|\\
		&\quad=\frac{1}{2M+1} \left| \sum_{j=-M}^M \mathbf 1_{A_{R,S}}(s-j)-\sum_{j=-M+1}^{M+1} \mathbf 1_{A_{R,S}}(s-j)\right|\\
		&\quad=\frac{|\mathbf 1_{A_{R,S}}(s+M)-\mathbf 1_{A_{R,S}}(s-M-1)|}{2M+1}\le \frac{1}{2M+1},\\
		&|w_{B_r(0),S_2}(\ell)-w_{B_r(0),S_2}(\ell-1)|\\
		&\quad=\frac{|\mathbf 1_{L}(\ell+M)-\mathbf 1_{L}(\ell-M-1)|}{2M+1}\le \frac{1}{2M+1}.
	\end{aligned}
\end{equation*}
Consequently, for all $s,\ell,n\in\mathbb Z$, we have 
\begin{equation*}
	\begin{aligned}
		|w_{R,S}(s)-w_{R,S}(s-n)|&\le \sum_{t=\min\{s,s-n\}+1}^{\max\{s,s-n\}} |w_{R,S}(t)-w_{R,S}(t-1)|\\
		&\le \frac{|n|}{2M+1},\\
		|w_{B_r(0),S_2}(\ell)-w_{B_r(0),S_2}(\ell-n)|&\le \sum_{t=\min\{\ell,\ell-n\}+1}^{\max\{\ell,\ell-n\}} |w_{B_r(0),S_2}(t)-w_{B_r(0),S_2}(t-1)|\\
		&\le \frac{|n|}{2M+1}.
	\end{aligned}
\end{equation*}
Since $|a-\sqrt{ab}|\le\sqrt{|a-b|}$ for $a,b\in [0,1]$, it follows that  
\begin{equation} \label{eq:sqrt-weight-estimate}
	\begin{aligned}
		\biggl| w_{R,S}(s)-\sqrt{w_{R,S}(s)w_{R,S}(s-n)} \biggr| &\le \sqrt{\frac{|n|}{2M+1}},\\
		\biggl| w_{B_r(0),S_2}(\ell)-\sqrt{w_{B_r(0),S_2}(\ell)w_{B_r(0),S_2}(\ell-n)} \biggr| &\le \sqrt{\frac{|n|}{2M+1}}.
	\end{aligned}
\end{equation}

For $(y,k)\in\mathbb T\times K$, let 
\[
\pi_{(y,k)}\colon C(\mathbb T\times K)\rtimes_{R_\alpha\times\sigma}\mathbb Z\to\mathcal{B}(\ell^2(\mathbb Z))
\]
be the regular representation given by 
\begin{equation} \label{eq:regular-representations}
	\pi_{(y,k)}(h)\delta_t = h(R_\alpha^t(y),\sigma^t(k))\delta_t,\qquad \pi_{(y,k)}(u_n)\delta_t = \delta_{t+n}
\end{equation}
for $h\in C(\mathbb T\times K)$ and $n,t\in\mathbb Z$, where $(\delta_t)_{t\in\mathbb Z}$ is the standard orthonormal basis of $\ell^2(\mathbb Z)$. For a finite subset $S\subseteq\mathbb Z$, let $P_S$ be the orthogonal projection onto $\ell^2(S)$, and identify $M_S=\mathcal B(\ell^2(S))$ with the corner $P_S\mathcal B(\ell^2(\mathbb Z))P_S$. 

Choose $y_R\in R$ for each $(R,S)\in\mathscr{R}_0\cup\mathscr{R}_1$ and $k_D\in D$ for each $D\in\mathcal Q$, and define $\psi_i\colon C(\mathbb T\times K)\rtimes_{R_\alpha\times\sigma}\mathbb Z\to F^{(i)}$ by
\begin{equation} \label{eq:psi-maps}
	\begin{aligned}
		\psi_i(a)&=\bigoplus_{(R,S)\in\mathscr R_i}\bigoplus_{D\in\mathcal Q} W_{R,S}^{1/2} P_S\pi_{(y_R,k_D)}(a)P_S W_{R,S}^{1/2},\qquad i=0,1,\\
		\psi_2(a)&=\bigoplus_{D\in\mathcal Q} W_{B_r(0),S_2}^{1/2} P_{S_2}\pi_{(0,k_D)}(a)P_{S_2} W_{B_r(0),S_2}^{1/2},
	\end{aligned}
\end{equation}
where $W_{R,S}$ and $W_{B_r(0),S_2}$ are the diagonal positive contractions defined in \eqref{eq:digonal-weights}. Each block map is completely positive, since it is a compression of a representation followed by multiplication on both sides by the square root of a positive contraction; see \cite[Example~1.5.2]{MR2391387}. At the unit, the block maps take values $W_{R,S}$ and $W_{B_r(0),S_2}$, respectively. Hence they are contractive. The norm of a finite direct sum is the maximum of the block norms, so each $\psi_{i}$ is a c.p.c.\ map.

For $h\in C(\mathbb T\times K)$, the same formulas give 
\begin{equation*}
	\begin{aligned}
		\psi_i(h) &= \bigoplus_{(R,S)\in\mathscr R_i}\bigoplus_{D\in\mathcal Q}\left( \sum_{s\in S} w_{R,S}(s) h(R_\alpha^s(y_R),\sigma^s(k_D)) e_{s,s}\right),\qquad i=0,1,\\
		\psi_2(h) &= \bigoplus_{D\in\mathcal Q}\left( \sum_{\ell\in S_2} w_{B_r(0),S_2}(\ell)h(R_\alpha^\ell(0),\sigma^\ell(k_D))\, e_{\ell,\ell} \right).
	\end{aligned}
\end{equation*}
In particular, 
\begin{equation} \label{eq:psi-diagonal}
	\psi_i(C(\mathbb T\times K))\subseteq D_F^{(i)},\qquad i=0,1,2.
\end{equation}

To relate the incoming maps to the partition of unity, put 
\begin{equation} \label{eq:widetilde-q_i}
	\widetilde{q}_i=\varphi^{(i)}(\psi_i(1)),\qquad i=0,1,2.
\end{equation}
By \eqref{eq:product-covariance}, \eqref{eq:color-maps}, and \eqref{eq:psi-maps}, we have 
\begin{equation} \label{eq:widetilde-q_i-expansion}
	\begin{aligned}
		\widetilde q_i &=\sum_{(R,S)\in\mathscr R_i}\sum_{D\in\mathcal Q}\sum_{s\in S} w_{R,S}(s) (\theta_s(f_R)\otimes \mathbf 1_{\sigma^s(D)}),\qquad i=0,1,\\
		\widetilde q_2 &=\sum_{D\in\mathcal Q}\sum_{\ell\in S_2} w_{B_r(0),S_2}(\ell) (\theta_\ell(\rho)\otimes \mathbf 1_{\sigma^\ell(D)}).
	\end{aligned}
\end{equation}
Using \eqref{eq:separated-weights} and \eqref{eq:weight-supports}, together with \eqref{eq:sum-form-q_i}, gives 
\begin{equation*}
	\begin{aligned}
		\widetilde q_i &= \frac{1}{2M+1}\sum_{j=-M}^M\sum_{(R,S)\in\mathscr R_i}\sum_{D\in\mathcal Q}\sum_{t\in A_{R,S}} u_{t+j}f_{R,D}u_{t+j}^*\\
		&=\frac{1}{2M+1}\sum_{j=-M}^M u_jq_iu_j^*,\qquad i=0,1,\\
		\widetilde q_2 &= \frac{1}{2M+1}\sum_{j=-M}^M\sum_{D\in\mathcal Q}\sum_{t\in L} u_{t+j}\rho_Du_{t+j}^* = \frac{1}{2M+1} \sum_{j=-M}^M u_jq_2u_j^*.
	\end{aligned}
\end{equation*}
Therefore,
\begin{equation} \label{eq:widetilde-q_i-sum1}
	\sum_{i=0}^2 \widetilde{q}_i = \frac{1}{2M+1}\sum_{j=-M}^M u_j(q_0+q_1+q_2)u_j^*=1.
\end{equation}
The functions $\widetilde{q}_i$ need not equal $q_i$, but their sum still recovers the unit exactly.

\medskip
\noindent\textbf{Approximation estimate.}
The following lemma estimates the approximation error on the dense $*$-algebra $C_c(\mathbb Z,C(\mathbb T\times K))$ of $C(\mathbb T\times K)\rtimes_{R_\alpha\times\sigma}\mathbb Z$.

\begin{lemma} \label{lem:compression-estimate}
	Let $E_0\subseteq \mathbb Z$ be a finite symmetric subset containing $0$, put $m=\max_{n\in E_0}|n|$, and fix an integer $M\ge 1$. Consider the corresponding product towers constructed as in \eqref{eq:product-towers}, with $[-M-m,M+m]\cap\mathbb Z\subseteq E$. Let $\varphi^{(i)}$, $\psi_i$, and $\widetilde q_i$ be as in \eqref{eq:color-maps}, \eqref{eq:psi-maps} and \eqref{eq:widetilde-q_i}, respectively, for $i=0,1,2$. Let $\eta>0$, and let
	\[
	p=\sum_{n\in E_0} h_n u_n \in C_c(\mathbb Z, C(\mathbb T\times K)).
	\]
	Suppose that each $h_n$ has oscillation at most $\eta$ on every product tower level, meaning that $|h_n(x)-h_n(y)|\le\eta$ for every $n\in E_0$ and all $x,y$ in the same product tower level. Then 
	\[
	\|\widetilde q_ip-\varphi^{(i)}(\psi_i(p))\|\le |E_0|\eta+\sqrt{\frac{m}{2M+1}}\sum_{n\in E_0}\|h_n\|,\qquad i=0,1,2.
	\]
\end{lemma}

\begin{proof}
	We first expand $\psi_i(p)$ in matrix units. For $(y,k)\in \mathbb T\times K$, $n\in E_0$, and $t\in \mathbb Z$, the regular representation \eqref{eq:regular-representations} gives
	\[
	\pi_{(y,k)}(h_n u_n)\delta_t = h_n(R_\alpha^{t+n}(y),\sigma^{t+n}(k))\delta_{t+n}.
	\]
	Thus nonzero entries occur only at the positions $(t+n,t)$, $t\in\mathbb Z$. Writing $s=t+n$, the $(s,s-n)$-entry is $h_n(R_\alpha^s(y),\sigma^s(k))$. The two-sided compressions in \eqref{eq:psi-maps} therefore give
	\begin{equation*}
		\begin{aligned}
			\psi_i(p)=&\bigoplus_{(R,S)\in\mathscr R_i}\bigoplus_{D\in\mathcal Q}\Biggl(\sum_{n\in E_0}\sum_{\substack{s\in S\\s-n\in S}}\sqrt{w_{R,S}(s)w_{R,S}(s-n)}\\
			&\hspace{65pt}\cdot h_n(R_\alpha^s(y_R),\sigma^s(k_D))e_{s,s-n}\Biggr),\qquad i=0,1,\\
			\psi_2(p)=&\bigoplus_{D\in\mathcal Q}\Biggl(\sum_{n\in E_0}\sum_{\substack{\ell\in S_2\\\ell-n\in S_2}}\sqrt{w_{B_r(0),S_2}(\ell)w_{B_r(0),S_2}(\ell-n)}\\
			&\hspace{35pt}\cdot h_n(R_\alpha^\ell(0),\sigma^\ell(k_D))e_{\ell,\ell-n}\Biggr).
		\end{aligned}
	\end{equation*}
	The restricted sums only involve matrix units inside the corresponding matrix algebras. By \eqref{eq:weight-columns}, an omitted index with $s-n\notin S$ has $w_{R,S}(s)=0$, and an omitted index with $\ell-n\notin S_2$ has $w_{B_r(0),S_2}(\ell)=0$.
	
	By the definition of $f_{R,D}$ and $\rho_D$ in \eqref{eq:f_{R,D}-rho}, the covariance relation in \eqref{eq:product-covariance}, and the definition of the single-tower maps in \eqref{eq:single-tower-castle-map}, we obtain
	\begin{equation*}
		\begin{aligned}
			\Phi_{R\times D,S,f_{R,D}}(e_{s,s-n})&=(\theta_s(f_R)\otimes \mathbf 1_{\sigma^s(D)})u_n,\\
			\Phi_{B_r(0)\times D,S_2,\rho_D}(e_{\ell,\ell-n})&=(\theta_\ell(\rho)\otimes \mathbf 1_{\sigma^\ell(D)})u_n.
		\end{aligned}
	\end{equation*}
	Applying $\varphi^{(i)}$ from \eqref{eq:color-maps} to the above expansions of $\psi_i(p)$ and using \eqref{eq:widetilde-q_i-expansion}, we obtain
	\[
	\widetilde q_i p-\varphi^{(i)}(\psi_i(p))=\sum_{n\in E_0} c_n^{(i)} u_n,\qquad i=0,1,2,
	\]
	where
	\begin{equation*}
		\begin{aligned}
		c_n^{(i)}&=\sum_{(R,S)\in\mathscr R_i}\sum_{D\in\mathcal Q}\sum_{s\in S}(\theta_s(f_R)\otimes 1_{\sigma^s(D)})\\
		&\qquad \cdot(w_{R,S}(s)h_n - \sqrt{w_{R,S}(s)w_{R,S}(s-n)}h_n(R_\alpha^s(y_R),\sigma^s(k_D))),\quad i=0,1,\\	
		c_n^{(2)}&=\sum_{D\in\mathcal Q}\sum_{\ell\in S_2} (\theta_\ell(\rho)\otimes 1_{\sigma^\ell(D)})\\
		&\qquad\cdot(w_{B_r(0),S_2}(\ell)h_n - \sqrt{w_{B_r(0),S_2}(\ell)w_{B_r(0),S_2}(\ell-n)}h_n(R_\alpha^\ell(0),\sigma^\ell(k_D))).
	\end{aligned}
	\end{equation*}
	Each term in the sum defining $c_n^{(i)}$ vanishes outside the corresponding product tower level, and the product tower levels within each family are pairwise disjoint. Hence, at each $(x,k)\in \mathbb T\times K$, at most one summand defining $c_n^{(i)}$ is nonzero. To estimate a nonzero summand at $(x,k)$, we write 
	\begin{equation*}
		\begin{aligned}
			& w_{R,S}(s)h_n(x,k) -\sqrt{w_{R,S}(s)w_{R,S}(s-n)} h_n(R_\alpha^s(y_R),\sigma^s(k_D)) \\
			&\quad=w_{R,S}(s)(h_n(x,k)-h_n(R_\alpha^s(y_R),\sigma^s(k_D)))\\
			&\qquad+(w_{R,S}(s)-\sqrt{w_{R,S}(s)w_{R,S}(s-n)})h_n(R_\alpha^s(y_R),\sigma^s(k_D)), \\
			& w_{B_r(0),S_2}(\ell)h_n(x,k) -\sqrt{w_{B_r(0),S_2}(\ell)w_{B_r(0),S_2}(\ell-n)}h_n(R_\alpha^\ell(0),\sigma^\ell(k_D)) \\
			&\quad=w_{B_r(0),S_2}(\ell)(h_n(x,k)-h_n(R_\alpha^\ell(0),\sigma^\ell(k_D))) \\
			&\qquad+(w_{B_r(0),S_2}(\ell)-\sqrt{w_{B_r(0),S_2}(\ell)w_{B_r(0),S_2}(\ell-n)})h_n(R_\alpha^\ell(0),\sigma^\ell(k_D)).
		\end{aligned}
	\end{equation*}
	In either case, $h_n$ is evaluated at two points of the same product tower level. Hence the oscillation assumption and the bounds $0\le w_{R,S},w_{B_r(0),S_2}\le 1$ give an absolute value of at most $\eta$ for the first term. By \eqref{eq:sqrt-weight-estimate}, the second term has absolute value at most 
	\[
	\sqrt{\frac{|n|}{2M+1}} \|h_n\|.
	\]
	Since $f_{R,D}$ and $\rho_D$ are positive contractions, their translates under the product action take values in $[0,1]$. It follows that
	\[
	\|c_n^{(i)}\|\le \eta+\sqrt{\frac{|n|}{2M+1}} \|h_n\|,\qquad i=0,1,2.
	\]
	Finally, since each $u_n$ is unitary and $|n|\le m$ for $n\in E_0$, we obtain
	\[
	\|\widetilde q_ip-\varphi^{(i)}(\psi_i(p))\|\le \sum_{n\in E_0}\|c_n^{(i)}\|\le |E_0|\eta+\sqrt{\frac{m}{2M+1}}\sum_{n\in E_0}\|h_n\|,\qquad i=0,1,2.
	\]
\end{proof}

\medskip
\noindent\textbf{Rectangular refinements.}
We will use the following elementary lemma to choose $\delta$ and a clopen partition of $K$ for which the oscillation assumption in Lemma~\ref{lem:compression-estimate} holds.

\begin{lemma} \label{lem:rectangular-refinement}
	Let $Y$ be a compact Hausdorff space, and let $K$ be a compact zero-dimensional Hausdorff space. For every finite open cover $\mathcal{U}$ of $Y\times K$, there exist a finite open cover $\mathcal{V}$ of $Y$ and a finite clopen partition $\mathcal{P}$ of $K$ such that
	\[
	\{V\times C : V\in\mathcal{V},\, C\in\mathcal{P}\}
	\]
	refines $\mathcal{U}$.
\end{lemma}

\begin{proof}
	Fix $y\in Y$. For each $k\in K$, choose an open neighborhood $O$ of $y$ and a clopen neighborhood $C$ of $k$ such that $O\times C$ is contained in some member of $\mathcal U$. By compactness of $K$, choose finitely many of these clopen neighborhoods covering $K$. Let $V_y$ be the intersection of the open neighborhoods of $y$ associated with the selected clopen sets, and let $\mathcal P_y$ be a finite clopen partition of $K$ refining the selected clopen cover. Then $V_y$ is an open neighborhood of $y$, and $V_y\times C$ is contained in some member of $\mathcal U$ for every $C\in\mathcal P_y$. By compactness of $Y$, there are $y_1,\dots,y_m\in Y$ such that 
	\[
	\mathcal V=\{V_{y_1},\dots,V_{y_m}\}
	\]
	covers $Y$. Let $\mathcal P$ be a finite clopen partition of $K$ refining each of $\mathcal P_{y_1},\ldots,\mathcal P_{y_m}$. Given $V_{y_j}\in \mathcal V$ and $C\in\mathcal P$, choose $C'\in\mathcal P_{y_j}$ with $C\subseteq C'$. Then $V_{y_j}\times C\subseteq V_{y_j}\times C'$, which is contained in some member of $\mathcal U$.
\end{proof}

\medskip
\begin{proof} [Proof of Theorem~\ref{thm:upper-bound-diagonal-dimension}]
	Recall that every minimal circle homeomorphism is topologically conjugate to an irrational rotation; see \cite[Proposition~11.1.4 and Theorem~11.2.7(1)]{MR1326374}. The product of this conjugacy with $\operatorname{id}_K$ conjugates $\tau\times\sigma$ to $R_\alpha\times \sigma$. By Proposition~\ref{prop:topological-conjugacy}, we may therefore assume that $\tau=R_\alpha$. Write 
	\[
	A=C(\mathbb T\times K)\rtimes_{R_\alpha\times \sigma} \mathbb Z.
	\]
	We verify Definition~\ref{def:diagonal-dimension} with $d=2$. 
	
	Let $\mathfrak{F}\subseteq A$ be a nonempty finite subset, and let $\varepsilon>0$. By density of $C_c(\mathbb Z,C(\mathbb T\times K))$, choose an element $p_a$ of this subalgebra for each $a\in\mathfrak{F}$ with 
	\[
	\|a-p_a\|<\frac{\varepsilon}{4}.
	\]
	Choose a finite symmetric subset $E_0\subseteq \mathbb Z$ containing $0$ such that every $p_a$ is supported in $E_0$, and write
	\[
	p_a=\sum_{n\in E_0} h_{a,n}u_n,\qquad a\in\mathfrak{F},
	\]
	extending the coefficients by zero if necessary. Put $m=\max_{n\in E_0}|n|$. Choose a number $\eta>0$ and an integer $M\ge 1$ such that 
	\begin{equation} \label{eq:chosen-M-eta-estimate}
		3|E_0|\eta<\frac{\varepsilon}{4}\quad\text{and}\quad 3\sqrt{\frac{m}{2M+1}} \max_{a\in\mathfrak{F}}\sum_{n\in E_0}\|h_{a,n}\|<\frac{\varepsilon}{4},
	\end{equation}
	and set 
	\[
	E=[-M-m,M+m]\cap \mathbb Z.
	\]
	Since $\{h_{a,n}\}_{a\in\mathfrak{F}, n\in E_0}$ is a finite family of continuous functions on the compact space $\mathbb T\times K$, there is a finite open cover $\mathcal U$ of $\mathbb T\times K$ such that each $h_{a,n}$ has oscillation less than $\eta$ on every member of $\mathcal U$. By Lemma~\ref{lem:rectangular-refinement}, there are a finite open cover $\mathcal V$ of $\mathbb T$ and a finite clopen partition $\mathcal P$ of $K$ such that 
	\[
	\{V\times C : V\in\mathcal{V},\, C\in\mathcal{P}\}
	\]
	refines $\mathcal U$. By the Lebesgue number lemma, choose $0<\delta<1/4$ such that every subset of $\mathbb T$ of diameter less than $\delta$ is contained in some member of $\mathcal V$. 
	
	Now apply the construction \eqref{eq:product-towers} with $E$, $\delta$ and $\mathcal P$, and let $\mathcal Q$ be the resulting clopen partition of $K$. Recall that $\mathcal Q$ refines $\sigma^{-s}(\mathcal P)$ for every $s\in H$, where $H$ is the union of the shapes of all towers involved. Each product tower level in \eqref{eq:product-towers} is the product of a circle level of diameter less than $\delta$ and a clopen set $\sigma^s(D)$, where $D\in\mathcal Q$ and $s\in H$ is the index of the level. By the choice of $\delta$, the circle level lies in some $V\in\mathcal V$, and by the choice of $\mathcal Q$, the set $\sigma^s(D)$ lies in some $C\in\mathcal P$. Hence each product tower level is contained in some $V\times C$ and thus in some member of $\mathcal U$. Consequently, each $h_{a,n}$ has oscillation less than $\eta$ on every product tower level.
	
	Let $\varphi^{(i)}$, $\psi_i$, and $\widetilde q_i$ be as in \eqref{eq:color-maps}, \eqref{eq:psi-maps} and \eqref{eq:widetilde-q_i}, respectively. By Lemma~\ref{lem:compression-estimate}, for each $a\in\mathfrak{F}$, 
	\[
	\|\widetilde q_i p_a-\varphi^{(i)}(\psi_i(p_a))\|\le |E_0|\eta+\sqrt{\frac{m}{2M+1}}\sum_{n\in E_0}\|h_{a,n}\|,\qquad i=0,1,2.
	\]
	For the finite-dimensional algebras $F^{(i)}$ and their diagonal subalgebras $D_F^{(i)}$ defined in \eqref{eq:finite-dimensional-algebras}, set  
	\[
	F=\bigoplus_{i=0}^2 F^{(i)},\qquad D_F=\bigoplus_{i=0}^2 D_F^{(i)},
	\]
	and
	\[
	\psi=\bigoplus_{i=0}^2 \psi_i\colon A\to F,\qquad \varphi=\sum_{i=0}^2 \varphi^{(i)}\colon F\to A.
	\]
	Then $D_F$ is a maximal abelian $*$-subalgebra of $F$. Since each $\psi_i$ is c.p.c.\ and the norm on $F$ is the maximum of the norms on its summands, $\psi$ is a c.p.c.\ map; moreover, \eqref{eq:psi-diagonal} gives $\psi(C(\mathbb T\times K))\subseteq D_F$. Since each $\varphi^{(i)}$ is c.p.c.\ order zero, $\varphi$ is completely positive. Let $v\in F$ be a matrix unit with respect to $D_F$, and let $v_i$ denote its image under the canonical projection $F\to F^{(i)}$. Since $v^*v$ is a minimal projection in $D_F$, exactly one $v_i$ is nonzero. Otherwise, $v^*v$ would split into two nonzero orthogonal projections in $D_F$, contradicting minimality. This nonzero $v_i$ is a matrix unit with respect to $D_F^{(i)}$, so \eqref{eq:matrix-unit-normalizer} yields 
	\[
	\varphi(v)=\varphi^{(i)}(v_i)\in \mathcal N_A(C(\mathbb T\times K)).
	\]
	Thus conditions (i), (iii), (iv), and (v) of Definition~\ref{def:diagonal-dimension} hold. It remains to verify the approximation condition (ii).
	
	By \eqref{eq:widetilde-q_i-sum1}, 
	\[
	\varphi\psi(1)=\sum_{i=0}^2 \varphi^{(i)}(\psi_i(1))=\sum_{i=0}^2 \widetilde{q}_i=1.
	\]
	Thus $\varphi\psi$ is unital and completely positive, and hence contractive. For $a\in\mathfrak{F}$, the above identity, Lemma~\ref{lem:compression-estimate}, and \eqref{eq:chosen-M-eta-estimate} give 
	\begin{equation*}
		\begin{aligned}
			\|\varphi\psi(p_a)-p_a\|&= \|\sum_{i=0}^2(\varphi^{(i)}(\psi_i(p_a))-\widetilde{q}_ip_a)\|\\
			&\le \sum_{i=0}^2 \|\varphi^{(i)}(\psi_i(p_a))-\widetilde{q}_ip_a\|\\
			&\le 3|E_0|\eta+3\sqrt{\frac{m}{2M+1}}\sum_{n\in E_0}\|h_{a,n}\|\\
			&<\frac{\varepsilon}{2}.
		\end{aligned}
	\end{equation*}
	Consequently, 
	\begin{equation*}
		\begin{aligned}
			\|\varphi\psi(a)-a\|&\le \|\varphi\psi(a-p_a)\|+\|\varphi\psi(p_a)-p_a\|+\|p_a-a\|\\
			&<\frac{\varepsilon}{4}+\frac{\varepsilon}{2}+\frac{\varepsilon}{4}=\varepsilon.
		\end{aligned}
	\end{equation*}
	This verifies condition (ii) of Definition~\ref{def:diagonal-dimension} and completes the proof.
\end{proof}

Taking $K$ to be a one-point space in Theorem~\ref{thm:upper-bound-diagonal-dimension}, we obtain the following exact computation for minimal circle homeomorphisms. 

\begin{cor} \label{cor:dimensions-minimal-circle-homeo}
	Let $\tau\colon \mathbb{T}\to \mathbb{T}$ be a minimal homeomorphism, and equip $\mathbb{T}$ with the $\mathbb{Z}$-action generated by $\tau$. Then
	\[
	\dim_{\mathrm{tow}}(\mathbb{T},\mathbb{Z}) = \dim_{\mathrm{ftow}}(\mathbb{T},\mathbb{Z}) = \dim_{\mathrm{diag}}(C(\mathbb{T}) \subseteq C(\mathbb{T}) \rtimes_\tau \mathbb{Z}) = 2.
	\]
\end{cor}

\begin{proof}
	By Theorem~\ref{thm:upper-bound-diagonal-dimension} with $K$ a one-point space, 
	\[
	\dim_{\mathrm{diag}}(C(\mathbb{T}) \subseteq C(\mathbb{T}) \rtimes_\tau \mathbb{Z}) \le 2.
	\] 
	For the lower bound, since $\mathbb{T}$ is path-connected and $\mathbb Z$ contains an element of infinite order, \cite[Theorem 1.1]{MR5032552} shows that the equivariant asymptotic dimension of this action with respect to the family of finite subgroups is at least $2$. The action is free by minimality, and the only finite subgroup of $\mathbb{Z}$ is trivial. Hence this dimension agrees with the amenability dimension $\dim_{\mathrm{am}}(\mathbb{T},\mathbb{Z})$; see \cite[Remark~4.14]{MR3606454}. Therefore, 
	\[
	\dim_{\mathrm{am}}(\mathbb{T},\mathbb{Z}) \ge 2.
	\]
	Since $\mathbb{T}$ is compact metrizable, the estimates in \cite[Theorem~5.14]{MR4167017} and \eqref{eq:fine-diagonal} yield
	\begin{equation*}
		\begin{aligned}
			2 \le \dim_{\mathrm{am}}(\mathbb{T},\mathbb{Z}) &\le \dim_{\mathrm{tow}}(\mathbb{T},\mathbb{Z}) \le \dim_{\mathrm{ftow}}(\mathbb{T},\mathbb{Z}) \\
			&\le \dim_{\mathrm{diag}}(C(\mathbb{T}) \subseteq C(\mathbb{T}) \rtimes_\tau \mathbb{Z}) \le 2.
		\end{aligned}
	\end{equation*}
\end{proof}

\begin{rem}
	By Denjoy's theorem \cite{zbMATH03008694}, every orientation-preserving $C^2$ diffeomorphism of the circle with irrational rotation number is topologically conjugate to an irrational rotation. Thus Corollary~\ref{cor:dimensions-minimal-circle-homeo} gives a positive answer to Question~\ref{ques:B} for these diffeomorphisms. 
	
	In the non-minimal case, the rotation towers pulled back along a semiconjugacy need not satisfy the oscillation assumption of Lemma~\ref{lem:compression-estimate}, so the argument above does not settle Question~\ref{ques:B}.
\end{rem}

%-------------------------------------------%

\section{Fine tower and diagonal dimensions of products with Cantor systems} \label{sec:fine tower and diagonal dimensions of products with Cantor systems}

In this section, we show that for the diagonal $\mathbb Z$-action on the product of a minimal circle system with a Cantor system, the fine tower dimension and the diagonal dimension are both equal to two. The upper bound follows from Theorem~\ref{thm:upper-bound-diagonal-dimension}. We prove the matching lower bound for fine tower dimension in a more general setting and show that the same lower bound holds for tower dimension when the action on the zero-dimensional factor has a periodic point.

The bases of open towers on a product space need not be rectangles. The following lemma gives rectangular bases with connected first factors while preserving the $E$-Lebesgueness and the number of colors.

\begin{lemma} \label{lem:rectangular-fine-towers}
	Let $Y$ be an infinite locally connected compact Hausdorff space, and let $K$ be a nonempty compact zero-dimensional Hausdorff space. Let $G$ be a countable discrete group acting on both $Y$ and $K$, and equip $Y\times K$ with the diagonal action. Fix a finite set $E\subseteq G$ and an integer $d\ge 0$. Suppose that $\{(B_i,S_i)\}_{i\in I}$ is a finite $E$-Lebesgue family of open towers whose supports have chromatic number at most $d+1$. Then there is a finite clopen partition $\mathcal{Q}$ of $K$ such that $Y\times K$ admits a finite $E$-Lebesgue family of open towers
	\[
	\{(V_\lambda\times C_\lambda,S_\lambda)\}_{\lambda\in\Lambda}
	\]
	such that: 
	\begin{enumerate}[label=(\roman*)]
		\item each $V_\lambda$ is a nonempty connected open subset of $Y$;
		\item each $C_\lambda$ belongs to $\mathcal{Q}$;
		\item for every $\lambda\in\Lambda$, there is $i\in I$ such that $V_\lambda\times C_\lambda\subseteq B_i$ with $S_\lambda=S_i$ and $(S_i)_E^\circ\neq \emptyset$;
		\item the tower supports have chromatic number at most $d+1$.
	\end{enumerate}
	In particular, if the original tower levels refine an open cover of $Y\times K$, then so do the refined levels.
\end{lemma}

\begin{proof}
	Fix a coloring $c\colon I\to\{0,\dots,d\}$ such that the supports of towers of the same color are pairwise disjoint. Since the family is $E$-Lebesgue, the levels $sB_i$, with $i\in I$ and $s\in (S_i)_E^\circ$, form a finite open cover of $Y\times K$. By the shrinking lemma for finite open covers of normal spaces, there are compact sets $F_{i,s}\subseteq sB_i$, with $i\in I$ and $s\in (S_i)^\circ_E$, that still cover $Y\times K$. For each $i\in I$, put 
	\[
	F_i=\bigcup_{s\in (S_i)^\circ_E}s^{-1}F_{i,s}\subseteq B_i.
	\]
	Then each $F_i$ is compact, and 
	\begin{equation} \label{eq:compact-interior-cover}
		Y\times K=\bigcup_{i\in I}\bigcup_{s\in (S_i)_E^\circ} sF_i.
	\end{equation}
	We may discard the indices for which $F_i=\emptyset$. By local connectedness of $Y$ and zero-dimensionality of $K$, every point of $F_i$ has a neighborhood $J\times D\subseteq B_i$, where $J$ is connected and open in $Y$ and $D$ is clopen in $K$. By compactness of $F_i$, there are finitely many such rectangles $J_{i,r}\times D_{i,r}$, $1\le r\le N_i$, such that 
	\[
	F_i\subseteq\bigcup_{r=1}^{N_i}(J_{i,r}\times D_{i,r})\subseteq B_i.
	\]
	
	Let $\mathcal{Q}$ be the finite clopen partition of $K$ generated by all the sets $D_{i,r}$ with $i\in I$ and $1\le r\le N_i$. For $i\in I$ and $C\in\mathcal{Q}$, set 
	\[
	O_{i,C}=\bigcup_{\{r:C\subseteq D_{i,r}\}}J_{i,r}.
	\]
	Since $O_{i,C}$ is a finite union of connected open sets, it has finitely many connected components, each of which is open. For each such component $V$, take $V\times C$ as a refined base. By construction, 
	\[
	V\times C\subseteq B_i.
	\]
	For fixed $i$, these bases are pairwise disjoint. They also cover $F_i$: if $(y,k)\in F_i$, choose $r$ with $(y,k)\in J_{i,r}\times D_{i,r}$ and let $C\in\mathcal Q$ contain $k$; then $C\subseteq D_{i,r}$ by the choice of $\mathcal Q$, so $y$ lies in a component of $O_{i,C}$. 
	
	For each refined base $V\times C$ obtained from $B_i$, take the tower $(V\times C,S_i)$. Each of its levels is contained in the corresponding level of $(B_i,S_i)$. The resulting family of towers is $E$-Lebesgue by \eqref{eq:compact-interior-cover}: if $x\in sF_i$ for some $i\in I$ and $s\in (S_i)^\circ_E$, then $s^{-1}x$ lies in one of the refined bases obtained from $B_i$.
	
	Finally, assign color $c(i)$ to each refined tower obtained from $(B_i,S_i)$. Refined towers from distinct original towers of the same color have disjoint supports because the original supports are disjoint. Refined towers from the same original tower also have disjoint supports: their bases are pairwise disjoint, and translates at different indices lie in distinct levels of the original tower. Thus the chromatic number does not increase. The containment of each refined level in its original level also preserves any prescribed open-cover refinement.
\end{proof}

We now use Lemma~\ref{lem:rectangular-fine-towers} to obtain lower bounds for the fine tower dimension and the tower dimension.

\begin{thm} \label{thm:product-fine-lower-bound}
	Let $Y$ be an infinite, connected, locally connected compact Hausdorff space, and let $\tau\colon Y\to Y$ be a minimal homeomorphism. Let $K$ be a nonempty compact zero-dimensional Hausdorff space, and let $\sigma\colon K\to K$ be a homeomorphism. Equip $Y\times K$ with the diagonal $\mathbb Z$-action generated by $\tau\times\sigma$. Then
	\begin{enumerate}[label=(\roman*)]
		\item $\dim_{\mathrm{ftow}}(Y\times K,\mathbb Z)\ge 2$;
		\item if $\sigma$ has a periodic point, then $\dim_{\mathrm{tow}}(Y\times K,\mathbb Z)\ge 2$.
	\end{enumerate}
\end{thm} 

\begin{proof}
	Since $Y$ is connected, the minimal homeomorphism $\tau$ is \textit{totally minimal}, meaning that $\tau^m$ is minimal for every integer $m\ge 1$; see \cite[II.5.4]{MR1249063}. Indeed, for every $m\ge 1$, $Y$ admits a finite partition into clopen $\tau^m$-invariant subsets on each of which $\tau^m$ acts minimally; see \cite[II.9.6.7]{MR1249063}. Connectedness forces this partition to be trivial.
	
	If $\sigma$ has a periodic point, fix a periodic orbit $O$ of period $q$. Otherwise, set $q=1$. In either case, put 
	\[
	E=\{-q,0,q\}.
	\]
	We show that the failure of either assertion yields a finite $E$-Lebesgue family of open towers with at most two colors and with bases of the form $V\times C$, where $V$ is a nonempty connected proper open subset of $Y$.
	
	For (i), suppose that $\dim_{\text{ftow}}(Y\times K,\mathbb Z)\le 1$. Choose a finite open cover $\mathcal W$ of $Y$ by proper open subsets. By the definition of fine tower dimension, there is a finite $E$-Lebesgue family of open towers with at most two colors whose levels refine $\{W\times K:W\in\mathcal W\}$. Apply Lemma~\ref{lem:rectangular-fine-towers} to this family. Its refined levels still refine the same cover. In particular, if the level of index $s$ of a refined tower with base $V\times C$ is contained in $W\times K$, then $V\subseteq \tau^{-s}(W)\subsetneq Y$. Thus each first factor of a refined base is a nonempty connected proper open subset of $Y$.
	
	For (ii), assume that $\sigma$ has a periodic point. Suppose that $\dim_{\mathrm{tow}}(Y\times K,\mathbb Z)\le 1$. Take a finite $E$-Lebesgue family of open towers with at most two colors and intersect its bases with the invariant subspace $Y\times O$, leaving the shapes unchanged. The resulting family is still $E$-Lebesgue on $Y\times O$ and has at most two colors. Apply Lemma~\ref{lem:rectangular-fine-towers} to this family. For each refined tower $(V\times C,S)$, choose $s\in S_E^\circ$. Since $q\in E$, we have $s+q\in S$. The levels indexed by $s$ and $s+q$ have the same nonempty second factor because $\sigma^q$ fixes $O$ pointwise. Their disjointness therefore implies $V\cap \tau^q(V)=\emptyset$, so $V\subsetneq Y$.
	
	In either case, we obtain a finite $E$-Lebesgue family of open towers
	\[
	\{(V_\lambda\times C_\lambda,S_\lambda)\}_{\lambda\in\Lambda}
	\]
	on $Y\times K_0$, where $K_0=K$ in the first case and $K_0=O$ in the second. Each $V_\lambda$ is a nonempty connected proper open subset of $Y$, and each $C_\lambda$ is a nonempty clopen subset of $K_0$. As observed after \cite[Definition~4.3]{MR4167017}, a tower can be reindexed by translating its base and its shape without changing its levels or its supports, and such reindexings preserve the $E$-Lebesgueness. Since $\tau$ and $\sigma$ are homeomorphisms, each translated base is still the product of a nonempty connected proper open subset of $Y$ and a nonempty clopen subset of $K_0$. We may therefore assume that $0$ is the least element of every $S_\lambda$.
	
	Fix a coloring $c\colon \Lambda\to\{0,1\}$ such that towers of the same color form a castle. We claim that for every $\lambda\in\Lambda$ there exist $\mu\in\Lambda$ and an integer $p\ge 1$ such that
	\begin{equation} \label{eq:recursive}
		V_\lambda\subseteq \tau^p(V_\mu).
	\end{equation}
	
	To prove the claim, fix $\lambda$ and a point $k\in C_\lambda$. Let $y\in V_\lambda$. By the $E$-Lebesgueness, there are $\mu\in\Lambda$ and $p\in(S_\mu)_E^\circ$ such that 
	\[
	(y,k)\in \tau^p(V_\mu)\times \sigma^p(C_\mu).
	\]
	Since $-q\in E$, we have $p-q\in S_\mu$. As $0$ is the least element of $S_\mu$, it follows that $p\ge q\ge 1$. The point $(y,k)$ also lies in the base level $V_\lambda\times C_\lambda$. Since $p\ge 1$, these two levels have distinct indices, namely $(\lambda,0)$ and $(\mu,p)$. Since distinct levels of the same color are disjoint, we conclude that $c(\mu)=1-c(\lambda)$. Set 
	\[
	\Gamma_\lambda(k)=\{(\mu,p): \mu\in \Lambda,\ p\in(S_\mu)_E^\circ,\ c(\mu)=1-c(\lambda),\ k\in\sigma^p(C_\mu)\}.
	\]
	As shown above, $p\ge 1$ for every $(\mu,p)\in\Gamma_\lambda(k)$, and 
	\[
	V_\lambda = \bigcup_{(\mu,p)\in\Gamma_\lambda(k)} (V_\lambda\cap \tau^p(V_\mu)).
	\]
	The sets in this finite union are relatively open in $V_\lambda$, and they are pairwise disjoint: a point $y$ in two of them would give a point $(y,k)$ lying in two distinct product levels of the same color, contradicting their disjointness. Since $V_\lambda$ is nonempty and connected, one of these sets equals $V_\lambda$. This proves the claim. 
	
	Now choose $\lambda_0\in\Lambda$. Applying the claim in \eqref{eq:recursive} repeatedly, we obtain $\lambda_{j}\in \Lambda$ and integers $p_j\ge 1$ for each $j\ge 0$ with 
	\[
	V_{\lambda_j}\subseteq \tau^{p_j}(V_{\lambda_{j+1}}).
	\]
	Since $\Lambda$ is finite, there are $0\le r<s$ with $\lambda_r=\lambda_s$. Combining the corresponding inclusions, we obtain 
	\[
	V_{\lambda_r}\subseteq \tau^{P}(V_{\lambda_r}),\qquad \text{where}\ P=\sum_{j=r}^{s-1} p_j \ge 1.
	\]
	Set $F=Y\setminus V_{\lambda_r}$. Then $F$ is a nonempty proper closed subset of $Y$, and the above inclusion implies  
	\[
	\tau^P(F)\subseteq F.
	\]
	Hence $(\tau^{nP}(F))_{n\ge 0}$ is a decreasing sequence of nonempty closed subsets of $Y$. By compactness of $Y$, the intersection 
	\[
	M=\bigcap_{n=0}^\infty \tau^{nP}(F)
	\]
	is a nonempty closed subset of $Y$. Moreover, 
	\[
	M\subseteq F\subsetneq Y\quad\text{and}\quad \tau^P(M)=\bigcap_{n=1}^\infty \tau^{nP}(F)=M.
	\]
	Thus $M$ is a nonempty proper closed $\tau^P$-invariant subset of $Y$, contradicting the minimality of $\tau^P$.
\end{proof}

\begin{rem}
	Taking $K$ to be a one-point space in the second assertion of Theorem~\ref{thm:product-fine-lower-bound} gives $\dim_{\mathrm{tow}}(Y,\mathbb Z)\ge 2$, since the unique point of $K$ is fixed by $\sigma$ and the product action is conjugate to the action generated by $\tau$. For $Y=\mathbb T$, this gives an alternative proof of the lower bound $\dim_{\mathrm{tow}}(\mathbb T,\mathbb Z)\ge 2$ in Corollary~\ref{cor:dimensions-minimal-circle-homeo} without invoking \cite[Theorem~1.1]{MR5032552}.
\end{rem}

We conclude with an application that separates the tower dimension from both the fine tower dimension and the diagonal dimension on a one-dimensional space. The estimates in \eqref{eq:free-cantor-factor-bounds} do not guarantee a separation between the two tower dimensions when $d=1$. The following corollary shows that, for products of minimal circle homeomorphisms with Cantor systems, such a separation occurs precisely when the Cantor action is free. 

\begin{cor} \label{cor:products-dimensions}
	Let $\tau\colon \mathbb T \to \mathbb T$ be a minimal homeomorphism, and let $\sigma\colon K \to K$ be a homeomorphism of a Cantor space $K$. Equip $\mathbb T\times K$ with the diagonal $\mathbb Z$-action generated by $\tau\times\sigma$. Then 
	\[
	\dim_{\mathrm{ftow}}(\mathbb T\times K,\mathbb Z) = \dim_{\mathrm{diag}}(C(\mathbb{T}\times K) \subseteq C(\mathbb{T}\times K) \rtimes_{\tau\times\sigma} \mathbb{Z}) = 2.
	\]
	Moreover,
	\[
	\dim_{\mathrm{tow}}(\mathbb{T} \times K, \mathbb{Z})=
	\begin{cases}
		1, & \text{if } \sigma \text{ generates a free } \mathbb{Z}\text{-action on } K,\\
		2, & \text{otherwise}.
	\end{cases}
	\]
\end{cor}

\begin{proof}
	By \eqref{eq:fine-diagonal}, Theorem~\ref{thm:upper-bound-diagonal-dimension}, and Theorem~\ref{thm:product-fine-lower-bound}(i), 
	\[
	2 \le \dim_{\mathrm{ftow}}(\mathbb{T} \times K, \mathbb{Z})
	\le \dim_{\mathrm{diag}}\big(C(\mathbb{T} \times K) \subseteq C(\mathbb{T} \times K) \rtimes_{\tau\times\sigma} \mathbb{Z}\big) \le 2.
	\]
	If $\sigma$ generates a free $\mathbb Z$-action on $K$, then the projection $\mathbb T\times K\to K$ is a $\mathbb Z$-equivariant continuous surjection onto a free Cantor system, so $\dim_{\mathrm{tow}}(\mathbb{T} \times K, \mathbb{Z}) = 1$ by \eqref{eq:tower-dimension-one}. Otherwise, $\sigma$ has a periodic point, so Theorem~\ref{thm:product-fine-lower-bound}(ii) and the equality established above give $2\le \dim_{\mathrm{tow}}(\mathbb T\times K,\mathbb Z)\le \dim_{\mathrm{ftow}}(\mathbb T\times K,\mathbb Z)=2$.
\end{proof}

%---------------------------------------------------------%
\bibliographystyle{plain}
\bibliography{references}
\end{document}